\documentclass[12pt]{amsart}
\usepackage{amsmath,mathrsfs}
\usepackage[english]{babel}

\usepackage[T1]{fontenc}

\newcommand{\N}{{\mathbb N}}

\newcommand{\C}{{\mathbb C}}
\newcommand{\R}{{\mathbb R}}

\newcommand{\de}{\partial}

\newcommand{\iu}{{\rm i}}

\numberwithin{equation}{section}
\newtheorem{theorem}{Theorem}[section]
\newtheorem{proposition}[theorem]{Proposition}
\newtheorem{lemma}[theorem]{Lemma}
\newtheorem{corollary}[theorem]{Corollary}
\newtheorem{definition}[theorem]{Definition}
\theoremstyle{definition}
\newtheorem{remark}[theorem]{Remark}

\newcommand{\brm}{\begin{remark}\rm}
\newcommand{\erm}{\end{remark}}
\newcommand{\brms}{\begin{remark}\rm}
\newcommand{\erms}{\end{remark}}
\newcommand{\bte}{\begin{theorem}}
\newcommand{\ete}{\end{theorem}}
\newcommand{\bpr}{\begin{proposition}}
\newcommand{\epr}{\end{proposition}}
\newcommand{\ble}{\begin{lemma}}
\newcommand{\ele}{\end{lemma}}
\newcommand{\beq}{\begin{equation}}
\newcommand{\eeq}{\end{equation}}
\newcommand{\bdm}{\begin{displaymath}}
\newcommand{\edm}{\end{displaymath}}
\numberwithin{equation}{section}

\newcommand{\bos}{\begin{remark}\rm}
\newcommand{\eos}{\end{remark}}

\newcommand{\ben}{\begin{enumerate}}
\newcommand{\een}{\end{enumerate}}

\newcommand{\be}{\begin{equation}}
\newcommand{\ee}{\end{equation}}

\title[On Schr\"odinger systems with nonlocal nonlinearities]{On a class of Schr\"odinger
systems with \\ local and nonlocal nonlinearities - PART1}

\author[H.\ Hajaiej]{Hichem Hajaiej}

\thanks{Ipeit (Institut preparatoire aux etudes d'ingenieur de Tunis)
 2, Rue Jawaher Lel Nahru -
1089 Montfleury - Tunis,
 Tunisie.
E-mail: {\em hichem.hajaiej@gmail.com}}

\address{Ipeit (Institut preparatoire aux etudes d'ingenieur de Tunis)
\newline\indent
 2, Rue Jawaher Lel Nahru -
1089 Montfleury - tunis
 Tunisie}
\email{hichem.hajaiej@gmail.com}

\thanks{ The  author
was partially supported by the Tunisian ARUB project : {\em Analyse
Math\'ematique et Applications 04/UR/15-02. }}

\begin{document}



\begin{abstract}
In this  first  part, we  study existence and  uniqueness  of solutions of a general nonlinear Schr\"odinger  system in the presence  of diamagnetic  field,  local and nonlocal nonlinearities.
\end{abstract}
\maketitle



\section{Introduction and main results}

\subsection{Introduction}

In this paper, we aim to study the following Cauchy problem of an $m$-coupled nonlinear Schr\"{o}dinger equations with electro-magnetic potentials, local and nonlocal nonlinearities
\begin{equation}
    \label{AVgh}
\begin{cases}
\iu \de_t \Phi_j = L_{A} \Phi_j + V(x) \Phi_j
- g_j(|x|, |\Phi_1|^2,\dots, |\Phi_m|^2) \Phi_j - \sum\limits_{i =1}^m W_{ij} * h(|\Phi_i|)  \frac{h'(|\Phi_j|)}{|\Phi_j|} \Phi_j, & \\
\noalign{\vskip1pt}
\Phi_j(0,x)=\Phi_j^0(x), & \\
\noalign{\vskip3pt}
1 \leq j \leq m,
\end{cases}
 \end{equation}
where, for all $1 \leq j \leq m$, $\Phi_j^0 : \R^N \rightarrow \C$,
$\Phi_j : \R^+ \times \R^N \rightarrow \C$, $h: \R^+ \rightarrow \R^+$ continuous and non-decreasing,
$V: \R^N \rightarrow \R$ and
$A: \R^N \rightarrow \R^N$ represent  electric and
magnetic potentials satisfying suitable assumptions
that will be stated in the following.
The magnetic operator  $L_{A}$ is  defined as
$$
L_{A} \phi:= \left( \frac{\nabla}{\iu}-A(x)   \right)^2 \phi = -\Delta \phi  - \frac2\iu A(x) \cdot \nabla \phi + |A(x)|^2 \phi -\frac1\iu \operatorname{div} A(x) \phi.$$
The magnetic field $B$ is $B= \nabla \times A$ in $\R^3$ and can be thought (and identified) in general dimension as a $2$-form $\mathcal{H}^{B}$ of coefficients $(\partial_i A_j - \partial_j A_i)$. We will keep using the notation $B=\nabla \times A$ in any dimension $N$. \\
In various relevant cases, it is possible to write \eqref{AVgh} in the following  vectorial form
\begin{equation*}
\begin{cases}
i \frac{\de \Phi}{\de t}= F'_A(\Phi) \\
\Phi(0,x)=\Phi^0 (x)
\end{cases}
 \end{equation*}
where  $\Phi^0=(\Phi^0_1,\dots,\Phi^0_m)$ and we have set
\begin{align*}
F_A(\Phi)& =\frac12 \sum_{j=1}^ m \int \left| \left( \frac{\nabla}{\iu}-A(x)   \right) \Phi_j  \right|^2 \, dx
+ \frac12 \int V(x) |\Phi|^2 \, dx- \int G(|x|, |\Phi_1|^2,\dots,|\Phi_m|^2) \, dx   \\ \nonumber
& -  \frac12 \sum_{i, j=1}^{m} \iint W_{ij}(|x-y|)  h(|\Phi_{i}(x)|)  h(|\Phi_{j}(y)|) \, dx dy,
\end{align*}
where  $W_{ij}=W_{ji}$ and such as $h$ satisfy suitable assumptions that will be stated in the following.
 Moreover, we observe that $G : (0, \infty) \times \R^m_{+} \rightarrow \R$ satisfy the following conditions
\begin{equation*}
\frac{\de G}{\de s_j}= g_j(|x|, s_1^2,\dots, s_m^2),
\end{equation*}
for every $j=1,\dots,m$. We look for a soliton or standing wave of \eqref{AVgh}, namely a solution of the form
$\Phi(t,x)=(\Phi_1(t,x),\dots,\Phi_m(t,x))$, where for $1 \leq j \leq m$,
$ \Phi_j(t,x)= e^{i \lambda_jt} u_j(x)$, $\lambda_j$ real numbers and $u_j: \R^N \rightarrow \C$.
Therefore, $\mathcal U = (u_1,\dots,u_m)$ is a solution of the following $m \times m$ elliptic problem:
\begin{equation}\label{elliptic}
\begin{cases}
L_{A} u_j + (V(x)-\lambda_j) u_j - g_j(|x|, |u_1|^2,\dots, |u_m|^2) u_j - \sum\limits_{i=1}^m W_{ij} * h(|u_i|)  \frac{h'(|u_j|)}{|u_j|} u_j=0 &\\
\noalign{\vskip4pt}
\text{for all $1 \leq j \leq m$}.
\end{cases}
 \end{equation}

\bigskip

\noindent
We point out that the general Schr\"{o}dinger system (\ref{AVgh}) we aim to study contains, as particular cases,  physically meaningful  situations as in Section 1.1 in \cite{ServSquassina}.

\subsection{Preliminaries and Notations}

\medskip

\noindent Since we want that the composite
functions $$x \mapsto G(|x|, u_1(x), \cdots, u_m(x))$$
are  measurable on $\R^N$ for every $u_1, \dots, u_m \in M(\R^N)$,
where $M(\R^N)$ is the set of measurable functions on $\R^N$,
we deal with the following $G$ of Carath\`eodory type:

\begin{definition}
A function $G: (0, \infty) \times \R^m \rightarrow \R$ is an m-Carath\`eodory function if
\begin{itemize}
\item [($1$)] $G(\cdot, s_1, \dots, s_m): (0, \infty) \rightarrow  \R$ is measurable on $(0, \infty) \setminus \Gamma $, where $\Gamma$ is a subset of $(0, \infty)$ having one dimension measure zero, for all $ s_1, \dots, s_m \geq 0  $,

\item [($2$)] For all $1 \leq n \leq m$, every $(m-1)$ tuple $s_i \geq 0$ and $r \in (0, \infty) \setminus \Gamma$, the function
\begin{eqnarray*}
\R & \rightarrow & \R \\
s_n & \mapsto & G(r, \dots, s_n, \dots)
\end{eqnarray*}
is continuous on $\R$.
\end{itemize}
\end{definition}

\medskip

\noindent Throughout this paper we denote by  $\mathcal{H}_{A}^1= \mathcal{H}_{A}^1(\R^N) =  (H^1_A(\R^N))^m$ where $H^1_A= H^1_A(\R^N)$ is the Hilbert space defined as the closure of $C_c^{\infty}(\R^N; \C)$ under the scalar product
$$ (u, v)_{{H}_{A}^1}= \Re \int (Du \cdot \overline{Dv}+ u \overline{v}) \, dx    $$
where $Du=(D_1 u, \dots, D_N u)$ and $D_j = i^{-1} \partial_j - A_j(x)$, with induced norm
$$\| u  \|^{2}_{{H}^1_{A}}= \int \left|\frac1i \nabla u- A(x) u \right|^2 \, dx + \int  |u|^2 \, dx < \infty.$$
Recall that the diamagnetic inequality
\begin{equation}\label{DiamIneq}
| \nabla | u| | \leq \left|  \left(  \frac{\nabla}{i} - A(x) \right)u \right|
\end{equation}
holds for every $u \in H^1_A(\R^N)$. The space $\mathcal{H}_{A}^1$ is equipped with the standard norm
$\|\Phi    \|^2_{\mathcal{H}^1_A}=\| (\frac{\nabla}{i}-A(x)) \Phi   \|^2_{\mathcal{L}^2}+ \| \Phi  \|^2_{\mathcal{L}^2}$ where $\mathcal{L}^p=\mathcal{L}^p(\R^N)= (L^p(\R^N))^m$ and the norm $\| \Phi  \|^2_{\mathcal{L}^p}= \sum_{i=1}^m \| \Phi_i \|^2_{L^p}$ for every $\Phi=(\Phi_1, \dots ,\Phi_m) \in \mathcal{L}^p$. We denote by $L^q_w(\R^N)$ with $q > 1$ the weak $L^q$-space (see   \cite{LiebLoss}) defined as the set of measurable functions $f$ equipped with the norm
$$\| f  \|_{q, w}= \sup_{D \subset \R^N, \ \mathcal{M}(D) < \infty}  (\mathcal{M}(D))^{-1/q'} \int_D |f(x)| \, dx < \infty, \ \ \frac1q+\frac1q'=1,$$
where $\mathcal{M}$ denotes the Lebesque measure on $\R^N$. The dual space of $\mathcal{H}^1_{A}$ is denoted by $\mathcal{H'}_{A}$. We denote $ \mathcal{H}^1= \mathcal{H}^1(\R^N)=(H^1(\R^N))^m $ equipped with the standard norm $\|\Phi    \|^2_{\mathcal{H}^1}=\| \nabla \Phi   \|^2_{\mathcal{L}^2}+ \| \Phi  \|^2_{\mathcal{L}^2}$  and    $\mathcal{H}^{-1}(\R^N)=(H^{-1}(\R^N))^m $. Clearly, by (\ref{DiamIneq}) the following Lemma holds:

\begin{lemma}\label{embeddH1A}
The space $\mathcal{H}_{A}^1$ is continuously embedded in $\mathcal{L}^p$ for all $p \in [2, 2^*]$ where $2^*=\frac{2N}{N-2}$ for $N \geq 3$ and there exists $C > 0$ independent on A such that
$$\| u \|_{\mathcal{L}^p} \leq C \| u  \|_{\mathcal{H}_{A}^1}  $$
Furthermore, $  (\mathcal{L}^p)'= \mathcal{L}^{p'} \subset \mathcal{H'}_{A}$ where $p'$ denotes the conjugate of $p$.
\end{lemma}

\noindent Recall that  by $C(\overline{I}, X)$ is the space of continuous functions $I \rightarrow X$  equipped with the unoform norm when $I$ is bounded.    By $\mathcal{D}(I, X)$,  we denote the space $C_c^{\infty}(I, X)$  of the $C^{\infty}$ functions $I \rightarrow X$ with compact support in $I$, equipped with the uniform norm of all derivatives on $I$. By $L^p(I, X)$ the Banach space of measurable function $I  \rightarrow X$ such that the norm
\begin{equation*}
\| u  \|_{L^p(I, X)}= \begin{cases}
\left( \int_I \|   u(t) \| _{X}^p \, dt  \right)^{1/p}, \quad\hbox{if $1 \leq p  < \infty$}, \\
{\rm Esssup} \  \| u(t) \| _{X}^p, \quad\hbox{if $ p  = \infty$}
\end{cases}
\end{equation*}
is finite. We denote by $W^{m, p}(I, X)$
 the Banach space of measurable functions $I \rightarrow X$ such that  $\frac{d^j u}{dt^j} \in L^p(I, X)$
 for every $j=1, \dots, m$ equipped with the norm
 $$
\| u \|_{W^{m, p}(I, X)} = \sum_{j=1}^m \| \frac{d^j u}{dt^j}  \|_{L^p(I, X)}.
 $$
$C^{m, \alpha}(\overline{I}, X)$ for $0 \leq \alpha \leq 1$ is the Banach  space of uniformly continuous and bounded functions with all their derivatives respect to t such that
$$
\| u \|_{C^{m, \alpha}(\overline{I}, X)}= \| u \|_{W^{m, \infty}(I, X)} + \sup \left\{ \frac{\left| \frac{d^m u(t)}{dt^m} - \frac{d^m u(s)}{dt^m}  \right|}{|t-s|^{\alpha}}; \ \ s, t \in I \right\} < \infty.
$$
Furthermore, we denote by $\mathcal{L}(X, Y)$ the Banach space of linear, continuous operators from the Banach spaces $X$  and $Y$ equipped with the norm topology.

\section{Local well-posedness}

\subsection{Assumptions on the magnetic potential}

We suppose that $A$ is a smooth function, namely $A \in C^{\infty}(\R^N, \R^N)$ and
there exist some constant $C_{\alpha} > 0$, $\alpha \in \N^n$ such that:
\begin{itemize}
\item [$(A)$] $\forall \alpha \in \N^n$, $|\alpha | \geq 1$,   $\sup_{x \in \R^N} | \partial^{\alpha}_x A | \leq C_{\alpha}$
\item  [$(B)$]    $\exists \varepsilon >0$, $\forall |\alpha  | \geq 1$, $\sup_{x \in \R^N} | \partial_x^{\alpha} B| \leq C_{\alpha} \left\langle x \right\rangle^{-1-\varepsilon}.$
\end{itemize}

\subsection{Assumptions on the external potentials}

We suppose that the external potentials satisfy:

\begin{itemize}

\item [$(V)$] $V \in L^p(\R^N, \R)$ or $ V: \R^N \rightarrow \R$, $V \in L^p +L^{\infty}$ for some $p \geq 1$, $p \geq N/2$;

\item [$(W)$]  $W_{ij}: \R^+ \rightarrow \R^+$ is an even, non-increasing function, $W_{ij}(|x|) \in L^q_w(\R^N)$ with $q > \max\{1, N/4 \}$ and $W_{ij}=W_{ji}$ for all $i, j= 1, \dots , m$.

\end{itemize}

\subsection{Assumptions on the local nonlinearities}

On the local non\-li\-nea\-ri\-ties, we assume that

\medskip

\noindent (g) For every $j=1, \dots, m$, the complex valued functions $f_j(x, \Phi)=g_j(|x|, |\Phi_1|^2,......., |\Phi_m|^2)\Phi_j$ are  measurable in $x \in \R^N$ and con\-ti\-nuous in $\Phi \in \C^m$     almost everywhere on $\R^N$. Assume that there exist constant $C$ and $\alpha \in [0, \frac{4}{N-2})$ ($\alpha \in [0, \infty)$ if $N=2$) such that
$$
| f_j(x, \Psi)- f_j(x,\Phi)| \leq C \left( |\Phi|^{\alpha} + |\Psi|^{\alpha}    \right) |\Psi - \Phi|,
\quad\hbox{for almost all $x \in \R^N$ and all $\Phi, \Psi \in \C^m$}.
$$
Observe that $f_j(x, \textbf{0})=0$.

\bigskip




\noindent (G) There exists $K>0$ such that, for all $r>0$ and $s_1,\dots,s_m  \geq 0$, we have
$$
0 \leq G(r, s_1,\dots, s_m) \leq K \Big(  \sum_{j=1}^{m} s_j + \sum_{j=1}^{m} s_j^{\frac{l_j +2}{2}}   \Big), \qquad  0 < l_j <  \frac{4}{N-2}.
$$

\begin{remark}
The local term in the energy functional $F_A$ is finite thanks to (G). Indeed,  we have
$$
\int G(|x|, |\Phi_1|^2,\dots,|\Phi_m|^2) \leq K \int | \Phi |^2 + K \sum_{j=1}^{m} \int | \Phi_j |^{l_j +2}.
$$
For every $j=1, \dots, m$, by the Gagliardo-Nirenberg and Diamagnetic inequalities, we have that
$$
\| \Phi_j   \|_{l_j +2} \leq c \| \Phi_j\|^{1- \sigma_j}_{L^2}   \| \nabla | \Phi_j| \|_{L^2}^{\sigma_j}  \leq c  \| \Phi_j \|^{1- \sigma_j}_{L^2}   \|  \Phi_j \|_{{\mathcal{H}^1_A}}^{\sigma_j}
 $$
where $\sigma_j= \frac{N l_j}{2(l_j +2)}$. Since  $\sigma_j$ must belong to $[0, 1)$, we find that $0 < l_j < \frac{4}{N-2}$.

\end{remark}






\subsection{Assumptions on the nonlocal nonlinearities}

On the nonlocal nonlinearity, we assume:
\vskip1pt
\noindent ($h$) $h: \R^+ \rightarrow \R^+$ is $C^1$ and non-decreasing, $h(0)=0$ and there exist $C,D,E>0$ such that
$$
h(s) \leq C s^\mu, \quad
|h'(s)| \leq D s^{\mu-1},\quad
|h''(s)| \leq E s^{\mu-2},
$$
for all $s\in\R^+$, where
\begin{equation}
    \label{restrh-localwp}
2 \leq \mu \leq  \frac{6q-1}{2q}- \frac{N-2}{N}.
\end{equation}
Notice that this inequality is nonempty due to the condition $q>N/4$.

\begin{remark}
If we set $H(s)= \frac{h'(s)}{s}$ for all $s \in \R^+$, there is a positive constant $C>0$ such that
\begin{align}
    \label{contr-H}
    |H(|z|)z-H(|w|)w| & \leq C(|z|^{\mu-2}+|w|^{\mu-2})|z-w|,\quad\text{for all $z,w\in\C$,} \\
    \label{contr-h}
    |h(|z|)-h(|w|)| & \leq C(|z|^{\mu-1}+|w|^{\mu-1})|z-w|,\quad\text{for all $z,w\in\C$.}
\end{align}
This easily follows in light of the growth conditions of the maps $\{s\mapsto h'(s),h''(s)\}$ contained in assumption $(h)$.
For the proof of inequality~\eqref{contr-H}, see for instance~\cite[inequality (2-1) of Lemma 2.1]{damasc}
applied with $A(\eta):=H(|\eta|)\eta:\R^2\to\R^2$. The only required condition is~\cite[condition (1-3)]{damasc}, which
is indeed fulfilled in view of the growths for $h'$ and $h''$ assumed in $(h)$.
\end{remark}


\begin{remark}\label{stimeW}
Since assumptions $(h)$ and $(W)$ hold true, by the Hardy-Littlewood-Sobolev inequality
for weak $L^q$ kernels (cf.~\cite[formula (7), p.107]{LiebLoss}), by the Gagliardo-Nirenberg
and  the diamagnetic inequality, one can prove that the nonlocal term involved in
the energy functional $F_A$ is always finite. Indeed, we have
\begin{align*}
 & \iint W_{ij}(|x-y|)  h(|\Phi_{i}(x)|)  h(|\Phi_{j}(y)|) \, dx dy \leq   C \| W_{ij} \|_{L^q_w} \| \Phi_i \|^{\mu}_{L^{\frac{2q \mu}{2q-1}}}  \| \Phi_j  \|^{\mu}_{L^{\frac{2q \mu}{2q-1}}} \\
\noalign{\vskip4pt}
 &  \leq   C \| W_{ij} \|_{L^q_w} \| \nabla |\Phi_i| \|_{L^2}^{N \mu \left( \frac12-\frac{(2q-1)}{2q \mu} \right)} \|  \Phi_i \|_{L^2}^{\mu \left[ 1- N\left( \frac12-\frac{(2q-1)}{2q \mu}\right) \right]} \| \nabla |\Phi_j| \|_{L^2}^{N \mu \left( \frac12-\frac{(2q-1)}{2q \mu} \right)} \|  \Phi_j \|_{L^2}^{\mu \left[ 1- N\left( \frac12-\frac{(2q-1)}{2q \mu}\right) \right]} \\
\noalign{\vskip4pt}
 & \leq C \| W_{ij} \|_{L^q_w} \left\| \left( \frac{\nabla}{i}- A(x)  \right) \Phi_i \right\|_{L^2}^{N \mu \left( \frac12-\frac{(2q-1)}{2q \mu} \right)} \|  \Phi_i \|_{L^2}^{\mu \left[ 1- N\left( \frac12-\frac{(2q-1)}{2q \mu}\right) \right]}  \\
\noalign{\vskip4pt}
    &   \times \left\| \left( \frac{\nabla}{i}- A(x)  \right) \Phi_j \right\|_{L^2}^{N \mu \left( \frac12-\frac{(2q-1)}{2q \mu} \right)} \|  \Phi_j \|_{L^2}^{\mu \left[ 1- N\left( \frac12-\frac{(2q-1)}{2q \mu}\right) \right]} \\
  & \leq  C \| W_{ij} \|_{L^q_w} \left\| \left( \frac{\nabla}{i}- A(x)  \right) \Phi \right\|_{\mathcal{L}^2}^{N \mu \left( \frac12-\frac{(2q-1)}{2q \mu} \right)} \|  \Phi_i \|_{L^2}^{\mu \left[ 1- N\left( \frac12-\frac{(2q-1)}{2q \mu}\right) \right]}  \\
\noalign{\vskip4pt}
    &   \times \left\| \left( \frac{\nabla}{i}- A(x)  \right) \Phi \right\|_{\mathcal{L}^2}^{N \mu \left( \frac12-\frac{(2q-1)}{2q \mu} \right)} \|  \Phi_j \|_{L^2}^{\mu \left[ 1- N\left( \frac12-\frac{(2q-1)}{2q \mu}\right) \right]} \\
\noalign{\vskip4pt}
&  \leq C \| W_{ij} \|_{L^q_w} m_{i}^{\frac{\mu}{2} \left[ 1- N\left( \frac12-\frac{(2q-1)}{2q \mu}\right) \right]}   m_{j}^{\frac{\mu}{2} \left[ 1- N\left( \frac12-\frac{(2q-1)}{2q \mu}\right) \right]} \\
\noalign{\vskip4pt}
&  \times  \left\| \left( \frac{\nabla}{i}- A(x)  \right)
\Phi \right\|_{\mathcal{L}^2}^{2N \mu \left( \frac12-\frac{(2q-1)}{2q \mu} \right)}   ,
\end{align*}
where $m_i= \| \Phi_i(t) \|^2_{L^2}$ and $m_j= \| \Phi_j(t) \|^2_{L^2}$ for all $i, j =1, \dots, m$. Observe that,
in order to have $ 2 \leq 2q \mu /(2q-1) \leq 2^*  $, $\mu$ must belong to the range $ 2 \leq \mu \leq 2^*(2q-1)/2q$,
which is compatible with the one in $(h)$, which is smaller.
\end{remark}

\medskip
\subsection{Local well-posedness}

In this section, we want to establish the local well posedness of the  Cauchy problem (\ref{AVgh})
\begin{equation*}
    \begin{cases}
\iu \de_t \Phi_j - L_A \Phi_j - V(x) \Phi_j
+ g_j(|x|, |\Phi_1|^2,\dots, |\Phi_m|^2) \Phi_j    + \sum\limits_{i=1}^m W_{ij} * h(|\Phi_i|)  \frac{h'(|\Phi_j|)}{|\Phi_j|} \Phi_j=0, & \\
\noalign{\vskip4pt}
\Phi_j(0,x)=\Phi_j^0(x), & \\
\noalign{\vskip4pt}
\text{for all $1 \leq j \leq m$},
\end{cases}
 \end{equation*}
or equivalently,
\begin{equation}
    \label{VgC}
\begin{cases}
\iu \de_t \Phi - L_A \Phi + \tilde{g}(\Phi)=0, & \\
\noalign{\vskip4pt}
\Phi(0,x)=\Phi^0(x), & \\
\noalign{\vskip4pt}
\end{cases}
 \end{equation}
where, for every $j=1,......,m$, $\tilde{g}_j(\Phi)= - \tilde{g}_{1,j}(\Phi)+
 \tilde{g}_{2,j}(\Phi) + \tilde{g}_{3,j}(\Phi)$ with $\tilde{g}_{1,j}(\Phi)=V(x) \Phi_j$, $ \tilde{g}_{2,j}(\Phi)=g_j(|x|, |\Phi_1|^2,\dots, |\Phi_m|^2) \Phi_j $ and   $ \tilde{g}_{3,j}(\Phi)=  \sum\limits_{i=1}^m W_{ij} * h(|\Phi_i|)  \frac{h'(|\Phi_j|)}{|\Phi_j|} \Phi_j$. Observe  that $L_A$ is a self-adjoint, $\geq 0$ operator on $\mathcal{L}^2$, $i L_A$ is  skew-adjoint and generates a group of isometries $\{T(t) \}_{t \in \R}$ where $T(t)= e^{- i t L_A}$ in $\mathcal{L}^2$. Furthermore, the following lemma holds

\begin{lemma}\label{JepsilonA}
If $\varepsilon > 0$ and $1 \leq p < \infty$, then $(I+ \varepsilon L_A)^{-1}$ is continuous  $\mathcal{L}^p \rightarrow \mathcal{L}^p$ and
$\| (I+ \varepsilon L_A)^{-1} \|_{\mathcal{L}(\mathcal{L}^p, \mathcal{L}^p)} \leq 1 $.

\end{lemma}

\begin{proof}
We adapt the proof of Lemma 9.1.3 in \cite{Cazenave}. In our case, we deal with  $\mathcal{L}^p$ and the operator $A$  is $L_A$.
\end{proof}

\begin{proposition}

\noindent  By  assumptions  $(V)$, $(W)$, $(g)$ and $(h)$ for suitable choice of $r_k$ and  $\rho_k$ (see Examples 1 to 4 in Section 4.2 in  \cite{Cazenave}), we have that   $\tilde{g}_k, k=1, 2,3$ satisfy the following conditions

\begin{enumerate}
    \item \label{4.2.1}
    $\tilde{g}_k \in C( \mathcal{H}_A^1, \mathcal{H}'_A)$ and there exists
    $\tilde{G}_k \in   C^1(\mathcal{H}_A^1, \R)$ such that $\tilde{g}_k=\tilde{G}'_k$;
\item
\label{4.2.2}
 there exist $r_k,\  \rho_k \in [2, \frac{2N}{N-2})$
($r_k, \rho_k \in [2, \infty)$ if $n=1$, $r_k,\rho_k \in [2, \infty)$ if $N=2$)
 such that   $\tilde{g}_k :\mathcal{H}_A^1 \rightarrow \mathcal{L}^{\rho_k'} \hookrightarrow \mathcal{H}'_A$;
\item
for every $M > 0$, there exists $C(M) < \infty$ such that
\begin{equation}\label{4.2.3}
\|\tilde{g}_k(\Psi)-  \tilde{g}_k(\Phi)   \|_{\mathcal{L}^{\rho_k'}}  \leq   C(M) \|\Psi -\Phi   \|_{\mathcal{L}^{r_k}},  \quad\hbox{for $k=1,2,3$},
\end{equation}
 for every $\Psi, \ \Phi \in \mathcal{H}_A^1$ such that $\| \Psi \|_{\mathcal{H}_A^1} + \| \Phi \|_{\mathcal{H}_A^1} \leq M$;
\item for every $\Phi \in \mathcal{H}_A^1$ and $j=1,....,m$,
\begin{equation}\label{4.2.4}
 \Im(\tilde{g}_{k, j}(\Phi) \overline{\Phi_j})=0 \quad\hbox{a.e. on $\R^N$}.
\end{equation}
\end{enumerate}
\end{proposition}

\begin{proof}
We start with $\tilde{g}_1$. Condition (\ref{4.2.1}) is satisfied since for every $j=1,\dots,m$,  $\tilde{g}_{1,j}$ satisfies it on $H^1_A$ and  $H^{-1}_A$ and  $\tilde{G}_1(\Phi)= \frac12 \int V(x) |\Phi|^2$ for all $\Phi \in \mathcal{H}^1_A$. (\ref{4.2.2}) and (\ref{4.2.3}) follow since they are satisfied by each component $\tilde{g}_{1,j}$ for every $j=1,\dots,m$ in the spaces $L^s$ with $s=\rho'_1, r_1$ where $r_1=\rho_1= \frac{2p}{p-1}$ and by the definition of the norm in the spaces $\mathcal{L}^s$, the diamagnetic and Sobolev inequalities.  Indeed,
\begin{align*}
\| \tilde{g}_{1,j} (\Psi) - \tilde{g}_{1,j} (\Phi)    \|_{L^{\rho'_1}}^{\rho'_1} & \leq   C  \int | V|^{\rho'_1} |\Psi_j-\Phi_j|^{\rho'_1}  \leq C
 \Big( \int | V|^{p} \Big)^{\frac{\rho'_1}{p}}     \Big( \int |\Psi_j-\Phi_j|^{\frac{p \rho'_1}{p-\rho'_1}}      \Big)^{\frac{p-\rho'_1}{p}}                   \\
 \noalign{\vskip4pt}
& \leq   C \|  V  \|_{L^p}^{\rho'_1} \| \Psi_j-\Phi_j  \|_{L^{\frac{p \rho'_1}{p-\rho'_1}}}^{\rho'_1}
\end{align*}
for  $\rho_1 \in [2, 2^*)$ and $r_1$ which satisfies
\begin{equation*}
2 \leq r_1=\frac{p \rho'_1}{p-\rho'_1} = \frac{p \rho_1}{p(\rho_1) -\rho_1} \leq 2^*.
\end{equation*}
The choice of $\rho_1= \frac{2p}{p-1}$   leads to $r_1=\rho_1 \in [2, 2^*)$ since $p \geq \frac{N}{2}$.
Finally, \eqref{4.2.4} follows since $V$ is real valued.
Condition~\eqref{4.2.1} is satisfied by $\tilde{g}_2$ since each component $\tilde{g}_{2,j}$, for every $j=1,....,m$, satisfies it on $H^1_A$ and  $H^{-1}_A$ and $\tilde{G}_2(\Phi)=  \int G(|x|, |\Phi_1|^2,\dots,|\Phi_m|^2) $. (\ref{4.2.2}) and (\ref{4.2.3}) follow easily from the local Lipschitz assumption in $(g)$, diamagnetic inequality, Sobolev embedding $H^1 \hookrightarrow L^{\alpha+2}$ and for $r_2= \rho_2= \alpha+2$. Indeed, by assumption (g),
\begin{align*}
\|\tilde{g}_{2,j}(\Psi)-\tilde{g}_{2,j}(\Phi)    \|_{L^{\rho_2'}}^{\rho_2'} &= \int \left| g_j(|x|, |\Psi_1|^2,\dots, |\Psi_m|^2) \Psi_j - g_j(|x|, |\Phi_1|^2,\dots, |\Phi_m|^2) \Phi_j    \right|^{\rho_2'}   \\
\noalign{\vskip4pt}
& \leq  C \Big(  \|  \Psi \|^{\alpha \rho'_2}_{\mathcal{L}^{\theta}} +  \|  \Phi \|^{\alpha \rho'_2}_{\mathcal{L}^{\theta}}  \Big)        \left\|\Psi -\Phi   \right\|_{\mathcal{L}^{r_2}}^{\rho'_2} \leq C(M)   \left\|\Psi -\Phi   \right\|_{\mathcal{L}^{r_2}}^{\rho'_2}.
\end{align*}
where $\rho_2, r_2 \in [2, 2^*)$ satisfy
\begin{equation*}
2 \leq \theta= \frac{\alpha \rho'_2 r_2}{r_2-\rho'_2}\leq 2^*.
\end{equation*}
The choice of $r_2= \rho_2= \alpha+2$   leads to $2 \leq \alpha +2 <  2^*$ which is compatible with the range of $\alpha$. Also (\ref{4.2.4}) is obvious since each $g_{2,j}$ is real valued. Finally, we deal with $\tilde{g}_3$.     (\ref{4.2.1}) holds since   for every $j=1,....,m$,  $\tilde{g}_{3,j}$ satisfies it on $H^1_A$ and  $H^{-1}_A$ and  $\tilde{G}_{3}(\Phi)= \frac12 \sum_{i, j=1}^{m} \iint W_{ij}(|x-y|)  h(|\Phi_{i}(x)|)  h(|\Phi_{j}(y)|) \, dx dy$. (\ref{4.2.2}) holds since, for each component, by H\"older inequality, assumption (h) and  Hardy-Littlewood-Sobolev inequality, we have
\begin{align*}
\|\tilde{g}_{3,j}(\Phi)   \|_{L^{\rho_3'}}^{\rho_3'} & \leq  C \sum\limits_{i=1}^m  \int  | W_{ij} * h(|\Phi_i|) |^{\rho_3'} \left| \frac{h'(|\Phi_j|)}{|\Phi_j|} \Phi_j \right|^{\rho_3'} \leq C \sum\limits_{i=1}^m \int  | W_{ij} * h(|\Phi_i|) |^{\rho_3'} \left| h'(|\Phi_j|) \right|^{\rho_3'} \\
\noalign{\vskip4pt}
& \leq  C \sum\limits_{i=1}^m \Big( \int  | W_{ij} * h(|\Phi_i|) |^{\frac{\rho_3' \rho_3}{\rho_3-\rho'_3}}\Big)^{\frac{\rho_3-\rho'_3}{\rho_3}}       \Big( \int   \left| h'(|\Phi_j|) \right|^{\rho_3}     \Big)^{\frac{\rho_3'}{\rho_3}} \\
\noalign{\vskip4pt}
& \leq  C \sum\limits_{i=1}^m \Big( \int  | W_{ij} * h(|\Phi_i|) |^{\frac{\rho_3 }{\rho_3-2}}\Big)^{\frac{\rho_3-\rho'_3}{\rho_3}}
  \Big( \int   |\Phi_j|^{\rho_3(\mu -1)}     \Big)^{\frac{\rho_3'}{\rho_3}}  \\
  \noalign{\vskip4pt}
  & \leq  C    \sum\limits_{i=1}^m   \| W_{ij} * h(|\Phi_i|)  \|^{   \frac{\rho_3 }{\rho_3-1}}_{L^{\frac{\rho_3 }{\rho_3-2}} }   \|  \Phi_j       \|^{\rho'_3(\mu -1)}_{L^{\rho_3(\mu -1)}}
\leq C \sum\limits_{i=1}^m     \| W_{i, j} \|_{L^q_w}^{\rho'_3} \| \Phi_i  \|_{L^{ \mu p}}^{\rho'_3 \mu}  \|  \Phi_j\|^{\rho'_3(\mu -1)}_{L^{\rho_3(\mu -1)}}
\end{align*}
where $p=  \frac{\rho_3 q}{2 \rho_3 q-2q-\rho_3}    $ for  $\rho_3 \in [2, 2^*)$  which satisfies
\begin{equation*}
\begin{cases}
2 \leq  \rho_3(\mu -1) \leq 2^*  \\
\frac{2}{\mu}  \leq   p=  \frac{\rho_3 q}{2 \rho_3 q-2q-\rho_3}  \leq \frac{2^*}{\mu}.
\end{cases}
\end{equation*}
Observe that the choice of $\rho_3= \frac{4q}{2q-1}$ in the above inequalities leads to the restriction
$$
2 \leq \mu < 1+ \frac{2^*(2q-1)}{4q},
$$
which is compatible with the range of $\mu$ in condition (h).
Condition \eqref{4.2.3}  is also satisfied  since,  for each component,  we have
\begin{align*}
\|\tilde{g}_{3,j}(\Psi)-\tilde{g}_{3,j}(\Phi)    \|_{L^{\rho_3'}}^{\rho_3'} & \leq  C \sum\limits_{i=1}^m  \int
\Big| W_{ij} * h(|\Psi_i|)   \frac{h'(|\Psi_j|)}{|\Psi_j|} \Psi_j - W_{ij} * h(|\Phi_i|)   \frac{h'(|\Phi_j|)}{|\Phi_j|} \Phi_j   \Big|^{\rho_3'}  \\
& \leq C  \sum\limits_{i=1}^m (I_i+J_i),
\end{align*}
where, for $i=1,\dots,m$ we have set
\begin{align*}
I_i &=\int  \Big| W_{ij} * h(|\Psi_i|) \Big|^{\rho_3'}     \Big| \frac{h'(|\Psi_j|)}{|\Psi_j|} \Psi_j -\frac{h'(|\Phi_j|)}{|\Phi_j|} \Phi_j       \Big|^{\rho'_3}    \\
\noalign{\vskip4pt}
J_i &= \int    \Big| W_{ij} * (h(|\Psi_i|)-h(|\Phi_i|)      )  \Big|^{\rho_3'}   \Big| h'(|\Phi_j|) \Big|^{\rho'_3}.
\end{align*}
By  virtue of H\"older inequality, condition~\eqref{contr-H} related to (h), Hardy-Littlewood-Sobolev and Sobolev  inequalities, we have that
\begin{align*}
 I_i &=   \int  \Big| W_{ij} * h(|\Psi_i|) \Big|^{\rho_3'}     \Big| H( |\Psi_j|) \Psi_j - H(|\Phi_j|)  \Phi_j       \Big|^{\rho'_3} \\
   \noalign{\vskip4pt}
   & \leq  \Big( \int  \Big| W_{ij} * h(|\Psi_i|) \Big|^{\frac{\rho_3}{\rho_3 -2}} \Big)^{\frac{\rho_3 -2}{\rho_3 -1}}    \Big( \int   \Big| H( |\Psi_j|) \Psi_j - H(|\Phi_j|)  \Phi_j       \Big|^{\rho_3} \Big)^{\frac{\rho'_3}{\rho_3}} \\
   \noalign{\vskip4pt}
   & \leq  C \| W_{i, j} \|_{L^q_w}^{\rho'_3} \| \Psi_i  \|_{L^{ \mu p}}^{\rho'_3 \mu}
    \Big( \int   \big(|\Psi_j|^{\rho_3(\mu-2)}+|\Phi_j|^{\rho_3(\mu-2)}\big)|\Psi_j-\Phi_j|^{\rho_3} \Big)^{\frac{\rho'_3}{\rho_3}} \\
    \noalign{\vskip4pt}
& \leq   C \| W_{i, j} \|_{L^q_w}^{\rho'_3} \| \Psi_i  \|_{L^{ \mu p}}^{\rho'_3 \mu}  \Big\{      \|  \Psi_j \|^{\rho'_3 (\mu -2)}_{L^{\frac{(\mu -2) r_3 \rho_3}{r_3 - \rho_3}}} +\|  \Phi_j \|^{\rho'_3 (\mu -2)}_{L^{\frac{(\mu -2) r_3 \rho_3}{r_3 - \rho_3}}}\Big\}     \| \Psi_j -\Phi_j \|^{\rho'_3}_{L^{r_3}}      \\
\noalign{\vskip4pt}
& \leq  C(M) \| \Psi_j -\Phi_j \|^{\rho'_3}_{L^{r_3}},
\end{align*}
where $\rho_3 \in [2, 2^*)$ and $r_3 \in [2, 2^*)$ satisfy the following conditions
\begin{equation}\label{SecondSystem}
\begin{cases}
2 \leq   \mu p           \leq  2^*  \\
2 \leq     \frac{(\mu -2) r_3 \rho_3}{r_3 - \rho_3}         \leq  2^* .
\end{cases}
\end{equation}
Observe that, by the choice of $\rho_3= \frac{4q}{2q-1}$, the above inequalities are satisfied by
\begin{equation}
\frac{4q 2^*  }{2^*(2q-1)-4q(\mu-2)}  \leq r_3 \leq \frac{4q}{6q-1-2q \mu}.
\end{equation}
The range of $\mu$ in assumption (h) ensures that $r_3 \in [2, 2^*)$.  Dealing with the second term in the above sum,
by means of condition~\eqref{contr-h} related to $(h)$, we have
\begin{align*}
        J_i &=  \int    \Big| W_{ij} * (h(|\Psi_i|)-h(|\Phi_i|)      )  \Big|^{\rho_3'}   \Big| h'(|\Phi_j|) \Big|^{\rho'_3}    \\
& \leq     \Big( \int    \Big| W_{ij} * (h(|\Psi_i|)-h(|\Phi_i|)      )  \Big|^{\frac{\rho_3}{\rho_3-2}}    \Big)^{\frac{\rho_3-2}{\rho_3-1}}      \Big( \int \Big| h'(|\Phi_j|) \Big|^{\rho_3}    \Big)^{\frac{\rho'_3}{\rho_3}}  \\
          \noalign{\vskip4pt}
 & \leq   \|  W_{ij} * (h(|\Psi_i|)-h(|\Phi_i|)   \|^{\rho'_3}_{L^\frac{\rho_3}{\rho_3-2}}    \|  h'(|\Phi_j|) \|_{L^{\rho_3}}^{\rho'_3}\\
 \noalign{\vskip4pt}
 & \leq  C \| W_{ij} \|_{L^q_w}^{\rho'_3}       \|   (h(|\Psi_i|)-h(|\Phi_i|)     \|^{\rho'_3} _{L^p }
 \|  \Phi_j    \|^{\rho'_3(\mu-1)}_{L^{\rho_3(\mu-1)}} \\
 \noalign{\vskip4pt}
 & \leq  C \| W_{ij} \|_{L^q_w}^{\rho'_3}   \Big(  \|  \Phi_i  \|_{L^{\frac{(\mu-1) p r_3}{(r_3-p)}}}^{\rho'_3(\mu-1)} +
   \|  \Psi_i   \|_{L^{\frac{(\mu-1) p r_3}{(r_3-p)}}}^{\rho'_3(\mu-1)}   \Big) \times  \\
& \times  \| \Psi_i-\Phi_i   \|_{L^{r_3}}^{\rho'_3}  \|  \Phi_j    \|^{\rho'_3(\mu-1)}_{L^{\rho_3(\mu-1)}}
  \leq  C(M) \| \Psi_i-\Phi_i   \|_{L^{r_3}}^{\rho'_3}
\end{align*}
where $p= \frac{\rho_3 q}{2 \rho_3 q-2q-\rho_3}      $     and $\rho_3 \in [2, 2^*)$ and $r_3 \in [2, 2^*)$ satisfy the following conditions
\begin{equation*}
\begin{cases}
2 \leq  \rho_3(\mu -1) \leq 2^*  \\
2 \leq  \frac{(\mu-1)p r_3}{(r_3-p)}          = \frac{(\mu -1) \rho_3 q r_3}{2 \rho_3 q r_3 -2q r_3 -\rho_3 r_3 - \rho_3 q} \leq 2^*.
\end{cases}
\end{equation*}
Observe again that, choosing   $\rho_3= \frac{4q}{2q-1}$, the second inequalities are satisfied by the same value of $r_3$ found for the second equations in (\ref{SecondSystem}). So
\begin{equation*}
\|\tilde{g}_{3}(\Psi)-\tilde{g}_{3}(\Phi)    \|_{\mathcal{L}^{\rho_3'}} \leq C(M)
  \|\Psi -\Phi   \|_{\mathcal{L}^{r_k}}
\end{equation*}
\eqref{4.2.4} is obvious since $W_{i,j}$, $h$ and $h'$ are real-valued.
\end{proof}

\noindent So  $\tilde{g}$  satisfy (\ref{4.2.1}) to (\ref{4.2.4}) and, in particular, assumptions like $(4.2.1)-(4.2.4)$ in Section $4.2$ or like $(4.6.3)-(4.6.6)$ in Section 4.6 in Cazenave \cite{Cazenave}. Consequently, by Remark 4.2.9 and Remark 4.2.13 or Remark 4.6.3 in \cite{Cazenave}, we are able to establish the local well-posedness  of the above Cauchy problem  in $\mathcal{H}_A^1$. Observe that  we assume the  ``a priori'' information that solutions are unique since   uniqueness is  proved by methods which are strictly related to the type of nonlinearity and,   in the following, we establish a result which ensures it. Now, we give the details of the proof. Recall that the energy functional is defined as
\begin{align*}
 F_A (\Phi) =&  \frac12 \sum_{j=1}^m \int \left| \left( \frac{\nabla}{\iu}-A(x)   \right) \Phi_j  \right|^2 \, dx
+ \frac12 \int V(x) | \Phi  |^2\, dx- \int G(|x|, |\Phi_1|^2,\dots,|\Phi_m|^2) \, dx \\
&  -
 \frac12 \sum_{i, j=1}^{m} \iint W_{ij}(|x-y|)  h(|\Phi_{i}(x)|)  h(|\Phi_{j}(y)|) \, dx dy \\
 =& \frac12 \left\| \left( \frac{\nabla}{i}-A(x)  \right) \Phi      \right\|^2_{\mathcal{L}^2}     + \tilde{G}_1(\Phi) -\tilde{G}_2(\Phi)-\tilde{G}_3(\Phi)
\end{align*}
and we denote by  $\tilde{G}(\Phi)= \tilde{G}_1(\Phi) -\tilde{G}_2(\Phi)- \tilde{G}_3(\Phi)$ for every $\Phi \in \mathcal{H}_A^1 $. It follows that $F_A \in C^1(\mathcal{H}_A^1, \R  )$ and that
$$F_A'(\Phi)= L_{A } \Phi- \tilde{g}(\Phi).   $$

\medskip

\noindent \textbf{Remark.}
In  assumption (\ref{4.2.1}), we require that $\tilde{g}: \mathcal{H}_A^1 \rightarrow \mathcal{H}'_A$ as $L_A$ does and that $\tilde{g}$ is the gradient of some functional $\tilde{G}$ since we can  define the energy. Indeed, the conservation of energy is essential in our proof of local existence.  (\ref{4.2.2}) requires that $\tilde{g}$ is slightly better than a mapping $\mathcal{H}_A^1  \rightarrow \mathcal{H}'_A $ and assumption (\ref{4.2.3}) is a form of local Lipschitz condition. Finally,  (\ref{4.2.4}) implies the conservation of charge that is essential for our proof.

\medskip

\begin{proposition}\label{solLinfinito}
Let  $A$ satisfy $(A)$ and $(B)$ with some constants $(C_{\alpha})_{\alpha \in \N^n}$ and assume (V), (W), (g) and (h) so that, in particular, $\tilde{g}_k, \  k=1,2,3$ satisfy assumptions (\ref{4.2.1})-(\ref{4.2.4}). For every $M > 0$, there exists $T(M)$ (depending only on $M$) and the $(C_{\alpha})$'s with the following property. For every $\Phi^0 \in \mathcal{H}_A^1$ such that $\| \Phi^0    \|_{\mathcal{H}_A^1} \leq M$, there exists a solution $\Phi \in L^{\infty} (I, \mathcal{H}_A^1) \cap W^{1, \infty} ( I, \mathcal{H}'_A)$ of the problem
\begin{equation}
    \label{Vg2}
\begin{cases}
 \iu \de_t \Phi - L_A \Phi + \tilde{g}(\Phi)=0, & \\
\noalign{\vskip4pt}
\Phi(0,x)=\Phi^0(x)
\end{cases}
 \end{equation}
with $I_M=(-T(M), T(M))$. In addition,
\begin{equation}\label{boundedLinfinito}
\| \Phi(t) \|_{L^{\infty} ((-T(M), T(M)); \mathcal{H}_A^1)  } \leq 2M.
\end{equation}
Furthermore,
\begin{equation}\label{consCharge2}
 \| \Phi(t) \|_{\mathcal{L}^{2}}=  \| \Phi^0 \|_{\mathcal{L}^{2}}
\end{equation}
\begin{equation}\label{consEnergy2}
F_A(\Phi(t)) \leq F_A(\Phi^0)
\end{equation}
for all $ t \in I_M=(-T(M), T(M))$.
\end{proposition}
\medskip

\noindent \textbf{Remark.} Note that both the equations  in (\ref{Vg2})  make sense respectively in $ \mathcal{H}'_A$ and  in $\mathcal{H}_A^1$. Indeed, since $\Phi \in W^{1, \infty} ((-T(M), T(M)) , \mathcal{H}'_A)$, by Remark 2.3.10 in \cite{Cazenave} easily adapted to our case, $\Phi$ is continuous $[-T(M), T(M)] \rightarrow  \mathcal{H}'_A$; and so, since $\Phi \in L^{\infty} ( (-T(M), T(M)), \mathcal{H}_A^1)$, it follows that $\Phi: [-T(M), T(M)] \rightarrow \mathcal{H}_A^1$ is weakly continuous. Furthermore, it follows from the duality inequality $\| \Phi  \|^2_{\mathcal{L}^2} \leq \|\Phi \|_{\mathcal{H}'_A} \   \| \Phi  \|_{\mathcal{H}_A^1} $ that $\Phi \in C( [-T(M), T(M)], \mathcal{L}^2 )$. Furthermore, also  (\ref{consCharge2}) and (\ref{consEnergy2}) make sense.

\medskip
\noindent In order to prove  the above proposition, we establish the following  two elementary lemmas.

\begin{lemma}\label{lemmaC0mezzi}
 Let $I \subset \R$ be an interval. Then, for every $\Phi  \in L^{\infty} (I, \mathcal{H}_A^1) \cap W^{1, \infty} ( I, \mathcal{H}'_A)$, we have
$$\| \Phi(t)- \Phi(s) \|_{\mathcal{L}^2} \leq C |t-s|^{1/2}, \quad\mbox{for all $s,t \in I $}    $$
where $C=\max \lbrace \| \Phi  \|_{L^{\infty} (I, \mathcal{H}^1_A}), \|\Phi'  \|_{L^{\infty} (I, \mathcal{H}'_A)} \rbrace $
\end{lemma}

\noindent \textbf{Proof.} The result follows from Remark 2.3.10 in \cite{Cazenave} applied to the case of $X= \mathcal{H}_A^1$ and $p=\infty$ and the duality  inequality $\| \Phi   \|^2_{\mathcal{L}^2}  \leq \|  \Phi\|_{\mathcal{H}'_A} \| \Phi    \|_{\mathcal{H}_A^1}$ (see Remark 2.3.8 (iii) in \cite{Cazenave} adapted to our case).

\medskip

\begin{lemma}\label{lemmagG}
 Let  $\tilde{g}_k$ for  $k=1, 2, 3$ satisfy assumptions (\ref{4.2.1})-(\ref{4.2.4}). Then, after possibly modifying the function $C(M)$ independent of $A$ (see Lemma 4.2.5, pag.55 in \cite{Cazenave} for  $A=0$ and one equation)
\begin{itemize}
\item [($i$)] for all $k=1, 2,3$, it holds $$\|\tilde{g}_k(\Psi) - \tilde{g}_k(\Phi)   \|_{\mathcal{L}^{\rho'_k}} \leq C(M) \|\Psi-\Phi   \|^{\alpha_k}_{\mathcal{L}^2}$$
for every $\Phi, \Psi \in \mathcal{H}_A^1$ such that $\| \Psi  \|_{ \mathcal{H}_A^1} +\|\Phi  \|_{ \mathcal{H}_A^1} \leq M$;
\medskip

 \item [($ii$)] for all $k=1, 2,3$, it holds $$|\tilde{G}_k(\Psi)-\tilde{G}_k(\Phi)
| \leq C(M) \|\Psi-\Phi   \|^{\beta_k}_{\mathcal{L}^2}$$
for every $\Phi, \Psi \in \mathcal{H}_A^1$ such that $\| \Psi  \|_{ \mathcal{H}_A^1} +\|\Phi  \|_{ \mathcal{H}_A^1} \leq M$,
\end{itemize}
 with $\alpha_k=1- N \big( \frac12-\frac{1}{r_k} \big)$ and $\beta_k=1- N \big( \frac12-\frac{1}{\rho_k}\big)$.
\end{lemma}

\begin{proof} (i) follows from $(\ref{4.2.3})$, the definition of the norm in $\mathcal{L}^{r_k}$,  Lemma \ref{embeddH1A} and Gagliardo-Nirenberg inequality  adapted in such case, namely  $\| \Psi-\Phi   \|_{\mathcal{L}^{r_k}}  \leq  C \|  \Psi-\Phi \|_{\mathcal{L}^2}^{1-\tilde{\alpha}_k } \| \Psi-\Phi   \|_{\mathcal{H}_A^1}^{\tilde{\alpha}_k}$ where $\tilde{\alpha}_k = N ( 1/2-1/r_k)$ (see~\cite[Theorem 2.3.7]{Cazenave}). (ii) follows from the identity
$$ \tilde{G}_k(\Psi)-\tilde{G}_k(\Phi) = \int_0^1 \frac{d}{ds} \tilde{G}_k (s\Psi +(1-s)\Phi) \, ds = \int_0^1 \left\langle \tilde{g}_k (s \Psi +(1-s)\Phi), \Psi -\Phi  \right\rangle_{\mathcal{H}'_A, \mathcal{H}_A^1}\, ds $$
hypothesis \eqref{4.2.2} and $\| \Psi-\Phi  \|_{\mathcal{L}^{\rho_k}}  \leq  C \|  \Psi-\Phi\|_{   \mathcal{L}^2       }^{1-\tilde{\beta}_k} \| \Psi-\Phi    \|_{\mathcal{H}^1_A         }^{\tilde{\beta}_k}$      where $\tilde{\beta}_k = N( 1/2-1/\rho_k)$.
The fact that the constant $C(M)$ is independent of the
magnetic field $A$ follows from the uniformity of the constant
in Lemma \ref{embeddH1A}.
\end{proof}

\bigskip

\noindent \textbf{Proof of Proposition \ref{solLinfinito}.} Following Cazenave \cite{Cazenave}, the proof proceeds in three steps. In the first step, we approximate $\tilde{g}$ by a family of ``regular'' nonlinearities. The importance of such regularization procedure is due to the necessity of obtaining the energy inequality  (\ref{consEnergy2}).         Dealing with the choice of the type of regularization which can be different for each nonlinearity we deal with, we  apply $(I+ \varepsilon L_A)^{-1}$ so that
the  proof applies to our different nonlinearities $\tilde{g}_k$, $k=1,2,3$ and  works as well when $\Omega=\R^N$. So  we  are able to construct approximate solutions. In the second step,  we obtain uniform estimates on such solutions, by using the conservation laws and, in particular,  the conservation of energy, in order to pass,  in the third step,   to the limit in the approximate equation. Observe that, even if there is a little bitter difficulty in the passage  to the limit in the nonlinearity,  we can recover the conservation of charge by the limiting problem thanks to the global regularization.

\bigskip

\noindent \textbf{Step 1: Construction of approximate solutions.}   From now on, we consider $\Phi^0 \in \mathcal{H}_A^1$ and we set $M=\| \Phi^0 \|_{\mathcal{H}_A^1} $. Given a positive integer $n \in \N$,  we define
$$ J_n^A= \left(I+ \frac1n L_A \right)^{-1}  $$

\noindent so that, for every $f \in \mathcal{H}'_A$, $ J_n^A f \in \mathcal{H}_A^1$ is the unique solution of the following system
$$
\Phi - \frac1n L_A \Phi=f \quad\hbox{in $\mathcal{H}'_A $}.
$$
From the self-adjointness of $L_A$, it is possible to deduce the following  main properties of the self-adjoint operator $ J_n^A$ (see Section 2.4 and Theorem 4.6.1 in Section 4.6 in \cite{Cazenave}), namely
\begin{equation}\label{4.2.10}
\| J_n^A \|_{\mathcal{L}(\mathcal{H}'_A, \mathcal{H}_A^1)} \leq \frac1n;
\end{equation}
\begin{equation}\label{4.2.11}
\| J_n^A \|_{\mathcal{L}(X,X)} \leq 1, \quad\hbox{whenever $X= \mathcal{H}_A^1,  \mathcal{L}^2,  \mathcal{H}'_A$};
\end{equation}

\begin{equation}\label{4.2.12}
J_n^A \Phi  \rightarrow_{n \rightarrow \infty} \Phi  \quad\hbox{in X, for every $\Phi \in X$, whenever $X= \mathcal{H}_A^1,  \mathcal{L}^2,  \mathcal{H}'_A$};
\end{equation}
\begin{equation}\label{4.2.13}
\hbox{if $\Phi^n$ is bounded in $X$, whenever $X= \mathcal{H}_A^1,  \mathcal{L}^2,  \mathcal{H}'_A$, then} \\ \nonumber
 J_n^A \Phi^n- \Phi^n \rightharpoonup 0 \quad\hbox{in $X$, as $n \rightarrow \infty$.}
\end{equation}

\noindent From Lemma \ref{JepsilonA}, we have that

\begin{equation}\label{4.2.14}
\| J_n^A \|_{\mathcal{L}(\mathcal{L}^p,\mathcal{L}^p)} \leq 1, \quad\hbox{for $ 1 \leq p < \infty$}
\end{equation}
We define, for every $\Phi \in \mathcal{L}^2$, the approximation of $\tilde{g} $ as
\begin{align*}
\tilde{g}_n(\Phi)& =  J_n^A (\tilde{g}(J_n^A \Phi))=J_n^A ((\tilde{g}_1 + \tilde{g}_2 +\tilde{g}_3)(J_n^A \Phi)) \\ \nonumber
 & =  J_n^A (\tilde{g}_1 (J_n^A \Phi)+ \tilde{g}_2 (J_n^A \Phi)+\tilde{g}_3(J_n^A \Phi))                     \\ \nonumber
 & =   J_n^A (\tilde{g}_1 (J_n^A \Phi)) + J_n^A ( \tilde{g}_2 (J_n^A \Phi)) +   J_n^A (\tilde{g}_3(J_n^A \Phi)) \\ \nonumber
 & =   \tilde{g}_{1, n} (J_n^A \Phi) + \tilde{g}_{2, n} (J_n^A \Phi) + \tilde{g}_{3, n} (J_n^A \Phi)                                       \end{align*}
and, for every $\Phi \in \mathcal{H}_A^1$, the approximation of $\tilde{G}$
\begin{equation}
\tilde{G}_n (\Phi)=\tilde{G}(J_n^A \Phi) =\tilde{G}_1(J_n^A \Phi) + \tilde{G}_2(J_n^A \Phi) + \tilde{G}_3(J_n^A \Phi) = \tilde{G}_{1, n}( \Phi) + \tilde{G}_{2, n}( \Phi) + \tilde{G}_{3, n}( \Phi).
\end{equation}
We observe that from (\ref{4.2.10})  the above definition make sense. Furthermore,  by (\ref{4.2.10}) and (\ref{4.2.3}), we have that $\tilde{g}_n$ is Lipschitz continuous on bounded sets of $\mathcal{L}^2$, and by  (\ref{4.2.10}) and  (\ref{4.2.1}) that $\tilde{G}_n \in C^1(\mathcal{H}_A^1, \R)$ and $\tilde{G}'_n=\tilde{g}_n$.  From (\ref{4.2.4}), it follows that

$$ (\tilde{g}_n(\Phi), i \Phi)_{\mathcal{L}^2}=(\tilde{g}(J_n^A \Phi), J_n^A \Phi )_{\mathcal{L}^2} =0, \quad\hbox{for every $\Phi \in \mathcal{L}^2$}.    $$
Therefore, 
there exists a sequence $(\Phi^n)_{n \in \N}$  of functions of $C(\R,\mathcal{H}_A^1) \cap C^1(\R, \mathcal{H}'_A)   $ such that
\begin{equation}\label{ApproximateProbl}
\begin{cases}
i \Phi^n_t - L_A \Phi^n + \tilde{g}_n (\Phi^n)=0, \\
\Phi^n (0)= \Phi^0.
\end{cases}
\end{equation}
Furthermore,
\begin{equation}\label{conservNormaPhim}
\| \Phi^n(t)  \|_{\mathcal{L}^2}  =\| \Phi^0 \|_{\mathcal{L}^2}
\end{equation}
and
\begin{align}\label{conservEnergyPhim}
F_A^n(\Phi^n(t) )& = \frac12 \left\| \left( \frac{\nabla}{i}-A(x)     \right) \Phi^n(t)   \right\|^2_{\mathcal{L}^2}+ \tilde{G}_n(\Phi^n(t)) \\
&=\frac12 \left\|  \left( \frac{\nabla}{i}-A(x)     \right)
\Phi^0  \right\|^2_{\mathcal{L}^2}+ \tilde{G}_n(\Phi^0)= F_A^n(\Phi^0 ) \notag
\end{align}
for all $t \in \R$.

\medskip

\noindent \textbf{Step 2. Estimates on the sequence $\Phi^n$.}  We denote by $C(M)$ various constants depending only on $M$. Remark again that the independence on $A$ follows from the uniformity of the constants involved in Lemma  \ref{embeddH1A} and consequently Lemmas \ref{lemmaC0mezzi} and \ref{lemmagG}. Let
\begin{equation}\label{thetam}
\theta_n = \sup \{ \tau >0:    \| \Phi^n(t)  \|_{\mathcal{H}^1_A} \leq 2M \quad\hbox{on $(\tau, \tau)$}   \}.
\end{equation}
Note that by (\ref{4.2.11}) and (\ref{4.2.14}),
\begin{equation}\label{UnifmA}
\tilde{g}_n \quad\hbox{verifies (\ref{4.2.2}) and (\ref{4.2.3}) uniformly in $n \in \N$ and $A$.  }
\end{equation}
Therefore, by (\ref{ApproximateProbl})
\begin{equation}\label{supPhit}
\sup_{n \in \N}   \| \Phi^n_t  \|_{L^{\infty}(-T(M), T(M); \mathcal{H}'_A)} \leq C(M).
\end{equation}
From (\ref{thetam}) and  (\ref{supPhit}), we can apply  Lemma \ref{lemmaC0mezzi} so

\begin{equation}\label{diffL2}
 \| \Phi^n(t) - \Phi^n(s) \|_{\mathcal{L}^2} \leq C(M) |t-s|^{1/2}, \quad\hbox{for all $s, t \in (- \theta_n, \theta_n)$.}
\end{equation}
Applying (\ref{conservNormaPhim}), (\ref{conservEnergyPhim}), Lemma \ref{lemmagG} (ii), (\ref{thetam}) and (\ref{diffL2}) for $s=0$, we obtain

\begin{align}\label{stimePhim}
\| \Phi^n(t)  \|_{\mathcal{H}^1_A}^2 &\leq  \| \Phi^0  \|_{\mathcal{L}^2}+ \left\| \left( \frac{\nabla}{i}-A(x)     \right) \Phi^0 \right\|_{\mathcal{L}^2} + 2 \left| \tilde{G}_n(\Phi^n(t))- \tilde{G}_n(\Phi^0) \right| \\ \nonumber
& \leq    \| \Phi^0  \|_{\mathcal{H}^1_A}^2 + 2 \sum_{k=1}^3 \left| \tilde{G}_{k,n}(\Phi^n(t))- \tilde{G}_{k, n}(\Phi^n(0)) \right|                                  \\ \nonumber
&  \leq  \| \Phi^0 \|_{\mathcal{H}^1_A}^2 + C(M) \sum_{k=1}^3 |t|^{\beta_k/2} \\
& \leq  \| \Phi^0 \|_{\mathcal{H}^1_A}^2 + C(M) |t|^{\beta/2}
\end{align}
where $\beta= \max \{ \beta_k: \ k=1, 2, 3 \}$ for all $ t \in (- \theta_n, \theta_n)$. If we define $T(M)$ by

$$C(M) T(M)^{\beta/2}=3 M^2,    $$
recalling that $\|\Phi^0  \|_{\mathcal{H}^1_A}$, it follows from (\ref{stimePhim}) that
$$ \| \Phi^n  \|_{L^{\infty}(-T, T; \mathcal{H}^1_A)} \leq 2M,  $$
for $T=\min\{ T(M), \theta_n  \}$. This implies that $T(M) \leq \theta_n$; and so

\begin{equation}\label{4.2.23}
\| \Phi^n \|_{L^{\infty}(-T(M), T(M); \mathcal{H}^1_A)} \leq 2M,
\end{equation}
and by (\ref{supPhit})
\begin{equation}\label{4.2.24}
\| \Phi^n_t  \|_{L^{\infty}(-T(M), T(M); \mathcal{H}'_A)} \leq C(M).
\end{equation}

\bigskip

\noindent \textbf{Step 3. Passage to the limit.} It follows from (\ref{4.2.23}) and (\ref{4.2.24}) and 
Proposition 2.3.13 (i)  
 in \cite{Cazenave} adapted to our case that there exists  $\Phi \in L^{\infty}(-T(M), T(M); \mathcal{H}^1_A) \cap W^{1, \infty}(-T(M), T(M); \mathcal{H}'_A)$ and a subsequence, which we still denote by $(\Phi^n)$ such that for all $t \in [-T(M), T(M)]$,

\begin{equation}\label{weakConvergPhim}
\Phi^n(t) \rightharpoonup \Phi(t) \quad\hbox{in $ \mathcal{H}^1_A$, as $n \rightarrow \infty$.}
\end{equation}
In addition, by (\ref{4.2.23}) and (\ref{weakConvergPhim}), Lemma \ref{lemmaC0mezzi}, (\ref{UnifmA}) and Lemma \ref{lemmagG}, we have that  $\tilde{g}_{k, n}(\Phi^n)$ is bounded in the space
$C^{0, \alpha_k /2}(-T(M), T(M); \mathcal{L}^{\rho'_k})$ for $k=1,2$ and $k=3$ for $\mu \in [2, 3)$ or in the space
$C^{0, \alpha_k (\mu-2) /2}(-T(M), T(M); \mathcal{L}^{\rho'_k})$ for $k=3$ and $\mu \geq 3$.
Therefore, it follows from Proposition 2.1.7 in \cite{Cazenave} adapted to our case that there exists $f_k$ which belongs to
$ C^{0, \alpha_k /2}(-T(M), T(M); \mathcal{L}^{\rho'_k})$ for $k=1,2$ and $k=3$ for $\mu \in [2, 3)$ or to $C^{0, \alpha_k (\mu-2) /2}(-T(M), T(M); \mathcal{L}^{\rho'_k})$ for $k=3$  and $\mu \geq 3$ and a subsequence, which we still denote by
$(\tilde{g}_{k, n}(\Phi^n))$ such that for all $t \in [-T(M), T(M)]$,
\begin{equation}\label{weakgm}
\tilde{g}_{k, n}(\Phi^n(t)) \rightharpoonup f_k(t) \quad\hbox{in $  \mathcal{L}^{\rho'_k}$, as $n \rightarrow \infty$.}
\end{equation}
On the other hand, it follows from (\ref{ApproximateProbl}) that for every $\Psi \in \mathcal{H}^1_A$ and for every $\phi \in \mathcal{D}(-T(M), T(M))$, we have
$$\int_{-T(M)}^{T(M)} \{ - \left\langle i \Phi^n, \Psi \right\rangle_{\mathcal{H}'_A, \mathcal{H}^1_A}  \phi'(t) -
 \left\langle L_A \Phi^n  - \sum_{k=1}^3 \tilde{g}_{k, n}(\Phi^n), \Psi         \right\rangle_{\mathcal{H}'_A, \mathcal{H}^1_A}  \phi(t)  \} \, dt=0 $$
Applying (\ref{weakConvergPhim}), (\ref{weakgm}) and the Dominated Convergence Theorem, it follows that
$$\int_{-T(M)}^{T(M)} \{ - \left\langle i \Phi, \Psi \right\rangle_{\mathcal{H}'_A, \mathcal{H}^1_A}  \phi'(t) - \left\langle L_A \Phi  - f, \Psi         \right\rangle_{\mathcal{H}'_A, \mathcal{H}^1_A}  \phi(t)  \} \, dt=0, $$
where $f=f_1+f_2+f_3$.  This means that $\Phi$ satisfies
\begin{equation}\label{LimitProbl}
\begin{cases}
i \Phi_t - L_A \Phi + f=0, \quad\hbox{for almost $t \in (-T(M), T(M))$}, \\
\Phi(0)= \Phi^0.
\end{cases}
\end{equation}
Now we prove  the following crucial result according to which the limit problem enjoys the conservation of charge.

\begin{lemma}
For all $t \in (-T(M), T(M))$, we have $\Im(f(t) \overline{\Phi(t)})=0$ almost everywhere on $\R^N$.
\end{lemma}

\begin{proof}
It's not so different  respect to the one in Lemma 4.2.6 in \cite{Cazenave}. Indeed, it's sufficient to show that for every bounded subsets $B$ of $\R^N$, we have for every $k=1, 2, 3,$

$$\left\langle f_k(t)|_{B}, i \Phi(t)|_B      \right\rangle_{ \mathcal{L}^{\rho'_k}(B),  \mathcal{L}^{\rho_k}(B)}=0.  $$
For simplicity, we omit the time dependence and we write
\begin{align*}
\left\langle f_k, i \Phi      \right\rangle_{ \mathcal{L}^{\rho'_k}(B),  \mathcal{L}^{\rho_k}(B)}&= \left\langle f_k-J_n^A (\tilde{g}_k(J_n^A \Phi^n)), i \Phi      \right\rangle +   \left\langle J_n^A (\tilde{g}_k(J_n^A \Phi^n)) -\tilde{g}_k(J_n^A \Phi^n) , i \Phi      \right\rangle    \\
& + \left\langle \tilde{g}_k(J_n^A \Phi^n), i (\Phi-\Phi^n) \right\rangle  +  \left\langle \tilde{g}_k(J_n^A \Phi^n), i (\Phi^n - J_n^A \Phi^n) \right\rangle \\
& + \left\langle \tilde{g}_k(J_n^A \Phi^n), i (J_n^A \Phi^n) \right\rangle \rightarrow_{n \rightarrow \infty} A+B+C+D+E.
\end{align*}
Note first that , by (\ref{weakgm}), $J_n^A (\tilde{g}_k(J_n^A \Phi^n))=\tilde{g}_{k, n}(\Phi^n) \rightharpoonup f_k$  in $ \mathcal{L}^{\rho'_k}$, hence in $\mathcal{L}^{\rho'_k}(B)$. Therefore, $A=0$. Next, we observe that $ \tilde{g}_k(J_n^A \Phi^n)$ is bounded in $\mathcal{L}^{\rho'_k}$. It follows from (\ref{4.2.13}) and (\ref{4.2.14}) that $J_n^A (\tilde{g}_k(J_n^A \Phi^n))-\tilde{g}_k(J_n^A \Phi^n)  \rightharpoonup 0   $  in $\mathcal{H}'_A$, hence in $\mathcal{L}^{\rho'_k}(B)$. Therefore, $B=0$. Since  by (\ref{weakConvergPhim}), $\Phi^n \rightharpoonup \Phi$ in $\mathcal{H}^1_A $, we have $\Phi^n \rightarrow \Phi$ in $\mathcal{L}^{\rho_k}(B)$. Since $\tilde{g}_k(J_n^A \Phi^n)$ is bounded in $\mathcal{L}^{\rho'_k}(B)$, it follows that $C=0$. By (\ref{4.2.11}) and (\ref{4.2.13}), $\Phi^n -J_n^A \Phi^n $ is bounded in $\mathcal{H}^1_A$ and converges weakly to $0$ in  $\mathcal{H}'_A$. It follows that $\Phi^n -J_n^A \Phi^n \rightarrow 0$ in $\mathcal{L}^{\rho_k}(B)$. Since $\tilde{g}_k(J_n^A \Phi^n)$ is bounded in  $\mathcal{L}^{\rho'_k}(B)$, it follows that $D=0$. Finally, $E=0$ by (\ref{4.2.4}). Hence the result.

\end{proof}

\bigskip

\noindent \textbf{End of the proof of Proposition \ref{solLinfinito}.}   Taking the $\mathcal{H}'_A - \mathcal{H}^1_A$ duality product of the first equation in (\ref{LimitProbl}) with $i \Phi$, it follows that

$$\frac{d}{dt} \| \Phi(t)  \|_{\mathcal{L}^2}=0, \quad\hbox{for all $t \in (-T(M), T(M))$;}   $$
and so
\begin{equation}\label{conservNormaL2}
 \| \Phi(t)  \|_{\mathcal{L}^2}=  \| \Phi^0  \|_{\mathcal{L}^2}.
 \end{equation}
It follows from (\ref{conservNormaPhim}), (\ref{conservNormaL2}) and Proposition 2.3.13 (ii) 
 in \cite{Cazenave}  adapted to our case that
\begin{equation}\label{convL2Phim}
 \Phi^n \rightarrow \Phi \quad\hbox{in $C(-T(M), T(M); \mathcal{L}^2)$.}
\end{equation}
Applying (\ref{4.2.23}), (\ref{convL2Phim}) and Gagliardo-Nirenberg inequality ( see Theorem 2.3.7  in \cite{Cazenave}), it follows that
\begin{equation}\label{convergLp}
 \Phi^n \rightarrow \Phi \quad\hbox{in $C(-T(M), T(M); \mathcal{L}^p)$, for every $2 \leq p < \frac{2N}{N-2} $.}
\end{equation}
It follows  from (\ref{4.2.3}), (\ref{4.2.12}) and (\ref{convergLp}) that
$$J_n^A (\tilde{g}_k(J_n^A \Phi^n))= \tilde{g}_{k,n} (\Phi^n(t)) \rightarrow \tilde{g}_k(\Phi(t)) \quad\hbox{in $\mathcal{L}^{\rho'_k}$, for all $ t \in  (-T(M), T(M))$.}  $$
Therefore, $f= \tilde{g}(\Phi)$ and so, $\Phi$ satisfies (\ref{Vg2}). (\ref{boundedLinfinito}) follows from (\ref{4.2.23}) and (\ref{consCharge2}) from (\ref{conservNormaL2}). It remains to prove (\ref{consEnergy2}). This follows from (\ref{conservEnergyPhim}), weak lower semicontinuity of the $\mathcal{H}^1_A$-norm and the fact that $\tilde{G}_n (\Phi^n(t)) \rightarrow \tilde{G}(\Phi(t))$ as $n \rightarrow \infty$ by (\ref{convergLp}) and Lemma \ref{lemmagG} (ii). This completes the proof.

\bigskip

\noindent Before proceeding further, we make the following definition.

\medskip

\begin{definition}
 In all what follows, we say that we have uniqueness for problem (\ref{Vg2}) if the following hold. For every interval $J$ containing $0$ and for every $\Phi^0 \in \mathcal{H}_A^1$, any two solutions of (\ref{Vg2}) in $ L^{\infty} (J, \mathcal{H}_A^1) \cap W^{1, \infty} ( J, \mathcal{H}'_A)$ coincide.
\end{definition}

\medskip

\noindent The main result of this section is the following.

\medskip

\begin{theorem}\label{solCC1}
 Let $A$  satisfies (A) and (B) and assume (V), (W), (g) and (h) so that, in particular, $\tilde{g}$ satisfy assumptions \eqref{4.2.1}-\eqref{4.2.4} and assume that we have uniqueness for problem~\eqref{Vg2}. Then the following properties hold.
\begin{enumerate}
 \item For every $\Phi^0 \in  \mathcal{H}_A^1(\R^N)$, there exists $T_{*}(\Phi^0), T^{*}(\Phi^0) > 0$ and there exists a unique, maximal solution $\Phi \in C(( - T_{*}(\Phi^0), T^{*}(\Phi^0)),\mathcal{H}_A^1) \cap C^1 (( - T_{*}(\Phi^0), T^{*}(\Phi^0)),\mathcal{H}'_A) $ of problem (\ref{Vg2}). $\Phi$ is maximal in the sense that if $ T^{*}(\Phi^0) < \infty $ (resp.,  $T_{*}(\Phi^0) < \infty)$, then $\| \Phi(t) \|_{\mathcal{H}^1_A} \rightarrow \infty$, as $t \uparrow T^{*}(\Phi^0)$ (resp., as $t \downarrow -T_{*}(\Phi^0) $);
\item in addition, we have conservation of charge and energy, that is
$$ \| \Phi(t) \|_{\mathcal{L}^{2}}=  \| \Phi^0 \|_{\mathcal{L}^{2}} \ \ \  F_A(\Phi(t)) = F_A(\Phi^0) $$
for all $t \in  ( -T_{*}(\Phi^0), T^{*}(\Phi^0))$;
\end{enumerate}

\end{theorem}

\noindent \textbf{Proof.} Following Cazenave \cite{Cazenave}, the proof proceeds in two steps. We first show that the solution $\Phi$ given by Proposition \ref{solLinfinito} belongs to $\Phi \in C(( - T_{*}(\Phi^0), T^{*}(\Phi^0)),\mathcal{H}_A^1) \cap C^1 (( - T_{*}(\Phi^0), T^{*}(\Phi^0)),\mathcal{H}'_A) $, and that we have conservation of energy. Next, we consider the maximality result.
\bigskip

\noindent \textbf{Step 1. Regularity.}  Let $I$ be an interval and let $\Phi \in L^{\infty} (I, \mathcal{H}_A^1) \cap W^{1, \infty} ( I, \mathcal{H}'_A)$ satisfy
$$i \Phi_t - L_A \Phi + \tilde{g}(\Phi)=0, \ \quad\mbox{for all $t \in I$} $$
We claim that $\Phi$ enjoys both conservation of charge and energy and that $\Phi \in C(I, \mathcal{H}_A^1) \cap C^1 ( I, \mathcal{H}'_A)$. To see this, consider
$$M=\sup\left\{ \| \Phi(t) \|_{\mathcal{H}_A^1}, t \in I      \right\},  $$
and let us first show that $\| \Phi(t) \|_{\mathcal{L}^2}$ and $F_A(\Phi(t))$ are constant on every interval $J \subset I$ of lenght at most $T(M)$, where $T(M)$ is given by Proposition \ref{solLinfinito}. Indeed, let $J$ be as above and let $\sigma, \tau \in J$. Let $\Phi^0= \Phi(\sigma)$ and let $\Psi$ be the solution of (\ref{Vg2}) given by Proposition \ref{solLinfinito}. $\Psi(\cdot-\sigma)$ is defined on $J$  and by uniqueness, $\Psi(\cdot-\sigma) = \Phi( \cdot)$ on $J$. By (\ref{consCharge2}) and (\ref{consEnergy2}), it follows in particular that
\begin{equation}\label{ConsSigmatau}
\| \Phi(\tau)  \|_{\mathcal{L}^2}=\| \Phi(\sigma)  \|_{\mathcal{L}^2}, \ \ \ \  F_A(\Phi(\tau))  \leq  F_A(\Phi(\sigma)).
\end{equation}
Now let $\Phi^0= \Phi(\tau)$ and let $Z$ be the solution of (\ref{Vg2}) given by Proposition \ref{solLinfinito}. $Z( \cdot -\tau)$ is defined on $J$ and by uniqueness, $Z( \cdot -\tau)= \Phi (\cdot)$ on $J$. By (\ref{consEnergy2}), it follows in particular that
$$F_A(\Phi(\sigma))  \leq  F_A(\Phi(\tau)).$$
Comparing with (\ref{ConsSigmatau}), this implies that both   $\| \Phi(t) \|_{\mathcal{L}^2}$ and $F_A(\Phi(t))$ are constant on  $J$. Since $J$ is arbitrary, it follows that
\begin{equation}\label{constantNorms}
\| \Phi(t) \|_{\mathcal{L}^2}=\| \Phi(s) \|_{\mathcal{L}^2} \quad\hbox{and}\ \ \  F_A(\Phi(t))=F_A(\Phi(s)), \quad\hbox{for all $s,t \in I$.}
\end{equation}
Furthermore, note that by Lemma \ref{lemmaC0mezzi}, $\Phi \in C^{0, 1/2} (\overline{I}, \mathcal{L}^2)$; and so, by Lemma \ref{lemmagG} (ii), the function $t \rightarrow \tilde{G}(\Phi(t))= \sum_{k=1}^3  \tilde{G}_k(\Phi(t)) $ is continuous $\overline{I} \rightarrow \R$. In view of (\ref{constantNorms}), it follows that
$\| \Phi(t) \|_{\mathcal{H}_A^1}$ is continuous $\overline{I} \rightarrow \R$. Therefore, by Lemma 2.1.5 in \cite{Cazenave} for $X=\mathcal{H}_A^1$, $\Phi \in C(\overline{I}, \mathcal{H}_A^1)$, and by the equation,  $\Phi \in C^1(\overline{I}, \mathcal{H}'_A)$.

\bigskip

\noindent \textbf{Step 2. Maximality.} Consider $\Phi^0 \in  \mathcal{H}_A^1$ and let
\begin{center}
\begin{equation*}
T^*(\Phi^0)= \sup \left\{ T > 0: \quad\hbox{there exists a solution of (\ref{Vg2}) on}\  [0, T]           \right\}
\end{equation*}
\begin{equation*}
 T_*(\Phi^0)= \sup \left\{ T > 0: \quad\hbox{there exists a solution of (\ref{Vg2}) on}\  [ -T, 0]           \right\}.
\end{equation*}
\end{center}

\medskip

\noindent By uniqueness and Step 1, there exists a solution
$$
\Phi \in C(( -T_{*}(\Phi^0), T^{*}(\Phi^0)),\mathcal{H}_A^1) \cap C^1 (( - T_{*}(\Phi^0), T^{*}(\Phi^0)),\mathcal{H}'_A)
$$
of \eqref{Vg2}. Suppose now that $T^{*}(\Phi^0) < \infty$ and assume that
there exists $M < \infty$ and a sequence $t_j \uparrow T^{*}(\Phi^0) $ such
that $\| \Phi(t_j) \|_{\mathcal{H}_A^1} \leq M$. Let $k$ be such that
$t_k +T(M) >  T^{*}(\Phi^0) $. By Proposition \ref{solLinfinito}    and Step 1 and starting
 from $\Phi(t_k)$, one can extend $\Phi$ up to $ t_k +T(M)  $, which is a contradiction with the maximality. Therefore, $\| \Phi(t) \|_{\mathcal{H}_A^1} \rightarrow \infty$, as $t \uparrow T^{*}(\Phi^0)$. One shows by the same argument that if $T_{*}(\Phi^0) < \infty$, then
 $\| \Phi(t) \|_{\mathcal{H}_A^1} \rightarrow \infty$, as $t \downarrow T_{*}(\Phi^0)$. Therefore, we have established statements (i) and (ii) of Theorem~\ref{solCC1}.

\bigskip

\noindent \textbf{Remark}. By Theorem \ref{solCC1}, under a priori uniqueness assumption, we have proved the well posedness of problem (\ref{Vg2}) in  $\mathcal{H}_A^1$, in particular, under assumptions
 (\ref{4.2.1}) through (\ref{4.2.4}) on  $\tilde{g}_k$  for $k=1, 2, 3$. We recall below a general sufficient condition for uniqueness by adapting Corollary  4.2.12    in   \cite{Cazenave}). It follows that
$$
\Psi(t)-\Phi(t)= i \int_0^t T(t-s) ( \tilde{g}(\Psi(s))-  \tilde{g}(\Phi(s)) \, ds,
$$
for all $t \in I$, where $T(t)$ is the propagator $e^{-it L_A}$. By assumptions (A) and (B)
on the potential and magnetic potentials, adapting the result in Yajima~\cite{Yajima} proved
for such $T(t)$, we have the following $\mathcal{L}^p$-$\mathcal{L}^q$ estimates

\begin{proposition}\label{stimePropagYajima}
Let $I_T=[-T, T] $. Then, for any $p$ such that $2 \leq p \leq \infty$ and $q$ conjugate to $p$, there exists a constant $C$ independent of $t$ such that for any $v \in \mathcal{L}^q$
\begin{equation}
    \label{decaytT}
\| T(t) v  \|_{\mathcal{L}^p} \leq \frac{C}{t^{N(\frac12-\frac1q)}} \| v \|_{\mathcal{L}^q}.
\end{equation}
\end{proposition}

\begin{corollary}
The conclusions of Theorem \ref{solCC1} holds true.
\end{corollary}

\begin{proof}
We have to prove that the uniqueness condition in $L^{\infty} (I, \mathcal{H}_A^1)
    \cap W^{1, \infty} ( I, \mathcal{H}'_A)$ is fulfilled.
    The argument follows the line of \cite[proof of Theorem 4.3.1]{Cazenave}.
 Let $I$ be an interval containing $0$
to be chosen sufficiently small. Let $\Psi, \Phi \in  L^{\infty} (I, \mathcal{H}_A^1)
\cap W^{1, \infty} ( I, \mathcal{H}'_A)$ be two solutions of equation~\eqref{Vg2}.
Let $r_i$ and $\rho_i$ the exponents for which the nonlinearity $g_i$ verifies the assumptions
of Theorem~\ref{solCC1}. Therefore, setting
$\frac{2}{q_j}=N(\frac12-\frac1{r_j})$ and $\frac{2}{\gamma_j}=N(\frac12-\frac1{\rho_j})$,
there exists $\delta>0$ such that
\begin{align*}
\|\Psi-\Phi\|_{L^{q_i}(I,L^{r_i})}\leq C\sum_{j=1}^m\|\tilde{g}(\Psi)- \tilde{g}(\Phi)\|_{L^{\gamma_j'}(I,L^{\rho_j'})}
\leq C(|I|+|I|^\delta)\sum_{j=1}^m\|\Psi- \Phi\|_{L^{q_j}(I,L^{r_j})},
\end{align*}
where the first inequality can be obtained by arguing
as in the proof of~\cite[Theorem 3.5.2(ii)]{Cazenave},
where the property in~\cite[Theorem 3.2.1]{Cazenave} is substituted by~\eqref{decaytT}, the estimate by Yajima.
In turn, adding the above inequality over $i=1,\dots,m$ and choosing the size of $|I|$ such that
$C(|I|+|I|^\delta)<1$ we get the inequality
$$
(1-C(|I|+|I|^\delta))\sum_{j=1}^m\|\Psi- \Phi\|_{L^{q_j}(I,L^{r_j})}\leq 0,
$$
yielding the desired conclusion.
\end{proof}

\medskip

\section{Global well-posedness}

\noindent We have established the local solvability of the Cauchy problem (\ref{Vg2}) in $\mathcal{H}^1_A$. In order to show that the solution $\Phi$ is global, namely that exists for all times, it is sufficient to establish a priori estimates on $\| \Phi(t) \|_{\mathcal{H}^1_A}$ by using the conservation laws (charge and energy) under some appropriate assumptions on the nonlinearities.

\medskip

\begin{theorem}\label{solGlobal}
Assume $(A)$, $(B)$, $(V)$, $(W)$, $(G)$ and $(h)$ with
$$
0 < l_j < \frac{4}{N},\,\,\,\qquad \mu < 2-\frac{1}{q}+\frac{2}{N},\qquad\inf_{x\in \R^N}V(x)>0.
$$
Let $\Phi^0 \in  \mathcal{H}_A^1(\R^N)$ be such that $\| \Phi^0 \|_{\mathcal{H}^1_A} \leq M$ and let
$$
\Phi \in C(( - T_{*}(\Phi^0), T^{*}(\Phi^0)),\mathcal{H}_A^1) \cap C^1 (( - T_{*}(\Phi^0), T^{*}(\Phi^0)),\mathcal{H}'_A)
$$
be the maximal solution of problem \eqref{Vg2} given by
Theorem~\ref{solCC1}. Then, $\Phi$ is global, namely
$ T_{*}(\Phi^0)= T^{*}(\Phi^0)= \infty$ and $\sup\{\|\Phi(t)\|_{\mathcal{H}^1_A}:t \in \R\}<\infty$.
\end{theorem}

\begin{proof}
Let $I_0= ( - T_{*}(\Phi^0), T^{*}(\Phi^0))$. By Theorem \ref{solCC1} (ii), we have the conservation of energy and charge, that is
$$ \| \Phi(t) \|_{\mathcal{L}^{2}}=  \| \Phi^0 \|_{\mathcal{L}^{2}} \ \ \  F_A(\Phi(t)) = F_A(\Phi^0) $$
for all $t \in I_0$. From the first equality we have that
$$ \| \Phi_j(t) \|_{\mathcal{L}^{2}}=  \| \Phi_j^0 \|_{\mathcal{L}^{2}} \quad\hbox{for all $j=1,....,m$}$$
and from the second
\begin{align*}
 F_A (\Phi(t)) &=  \frac12 \sum_{j=1}^m \int \left| \left( \frac{\nabla}{\iu}-A(x)   \right) \Phi_j(t)  \right|^2
+ \frac12 \int V(x) | \Phi(t)  |^2 \\
&- \int G(|x|, |\Phi_1(t)|^2,\dots,|\Phi_m(t)|^2)  \\
\noalign{\vskip4pt}
&  -
 \frac12 \sum_{i, j=1}^{m} \iint W_{ij}(|x-y|)  h(|\Phi_{i}(t)|)  h(|\Phi_{j}(t)|) \, dx dy  = F_A(\Phi^0) =C_0.
\end{align*}
Since $V$ is bounded from below we have that
\begin{align*}
C\| \Phi(t) \|^2_{\mathcal{H}^1_A} & \leq  \left\| \left( \frac{\nabla}{i}-A(x)   \right) \Phi(t)   \right\|^2_{\mathcal{L}^2} + \int V(x) |\Phi(t)|^2 \\
\noalign{\vskip4pt}
& \leq  C_0+ 2 \int G(|x|, |\Phi_1(t)|^2,\dots,|\Phi_m(t)|^2)  \\
& +\sum_{i, j=1}^{m} \iint W_{ij}(|x-y|)  h(|\Phi_{i}(t)|)  h(|\Phi_{j}(t)|) \, dx dy
\end{align*}
By assumptions $(G0)$-$(G1)$, we have that
\begin{align*}
\int G(|x|, |\Phi_1(t)|^2,\dots,|\Phi_m (t)|^2) & \leq K \int |\Phi(t)|^2 + K \int \sum_{j=1}^{m} | \Phi_j(t)|^{l_j +2} \\
& = K \| \Phi(t)  \|^2_{\mathcal{L}^2} + K  \sum_{j=1}^{m} \int  | \Phi_j(t)|^{l_j +2}
\end{align*}
For $j=1,...,m$, by the Gagliardo-Nirenberg inequality we have that:
$$
\| \Phi_j (t)  \|_{l_j +2} \leq c \| \Phi_j(t) \|^{1- \sigma_j}_{L^2}   \| \nabla | \Phi_j(t)| \|_{L^2}^{\sigma_j}, \ \ \sigma_j= \frac{N l_j}{2(l_j +2)}.
 $$
Now let $p_j= \frac{4}{N l_j}$ and $q_j$ is such that $\frac{1}{p_j}+ \frac{1}{q_j}=1$. Applying Young  and Diamagnetic Inequalities, we obtain
\begin{align*}
\| \Phi_j (t)  \|_{l_j +2}^{l_j +2} & \leq  C^{l_j +2}  \| \Phi_j(t) \|^{(1- \sigma_j)(l_j +2)}_{L^2}   \| \nabla | \Phi_j(t)| \|_{L^2}^{\sigma_j(l_j +2)}  \\
\noalign{\vskip4pt}
& \leq  \frac{1}{q_j} \{ \frac{C^{l_j +2}}{\varepsilon}  \| \Phi_j(t) \|^{(1- \sigma_j)(l_j +2)}_{L^2}   \}^{q_j} + \frac{N l_j}{4} \varepsilon^{\frac{4}{N l_j}} \| \nabla | \Phi_j(t)| \|_{L^2}^2 \\
\noalign{\vskip4pt}
& \leq  \frac{1}{q_j} \{ \frac{C^{l_j +2}}{\varepsilon}  m_{j}^{\frac{(1- \sigma_j)(l_j +2)}{2}}  \}^{q_j} + \frac{N l_j}{4} \varepsilon^{\frac{4}{N l_j}} \left\| \left( \frac{\nabla}{i}- A(x)   \right) \Phi_j(t) \right\|_{L^2}^2 \\
\noalign{\vskip4pt}
& \leq  \frac{1}{q_j} \{ \frac{C^{l_j +2}}{\varepsilon}       m_{j}^{\frac{(1- \sigma_j)(l_j +2)}{2}}  \}^{q_j} + \frac{N l_j}{4} \varepsilon^{\frac{4}{N l_j}} \left\| \left( \frac{\nabla}{i}- A(x)   \right) \Phi_j(t) \right\|_{\mathcal{L}^2}^2 \\
\noalign{\vskip4pt}
& \leq  \frac{1}{q_j} \{ \frac{C^{l_j +2}}{\varepsilon}  m_{j}^{\frac{(1- \sigma_j)(l_j +2)}{2}}   \}^{q_j} + \frac{N l_j}{4} \varepsilon^{\frac{4}{N l_j}} \left\|  \Phi_j(t) \right\|_{\mathcal{H}^1_A}^2.
\end{align*}
Consequently,
\begin{align*}
\int G(|x|, |\Phi_1 (t)|^2,\dots,|\Phi_m (t)|^2)  & \leq  K \| \Phi(t)  \|^2_{\mathcal{L}^2} + K \sum_{j=1}^{m}
 \frac{1}{q_j} \{ \frac{C^{l_j +2}}{\varepsilon}  m_{j}^{\frac{(1- \sigma_j)(l_j +2)}{2}}   \}^{q_j} \\
& + K
\left(  \sum_{j=1}^{m} \frac{N l_j}{4} \varepsilon^{\frac{4}{N l_j}} \right)  \left\|  \Phi(t)\right\|_{\mathcal{H}^1_A}^2 \\
\noalign{\vskip4pt}
& =   C_1+ C_2(\varepsilon) +  K \left(  \sum_{j=1}^{m} \frac{N l_j}{4} \varepsilon^{\frac{4}{N l_j}} \right)            \left\|  \Phi(t) \right\|_{\mathcal{H}^1_A}^2.
\end{align*}
Following the calculations done in Remark \ref{stimeW},
\begin{align*}
&   \sum_{i, j=1}^{m} \iint W_{ij}(|x-y|)  h(|\Phi_{i}(t)|)  h(|\Phi_{j}(t)|) \, dx dy  \\
\noalign{\vskip4pt}
 & \leq    \sum_{i, j=1}^{m}  \| W_{ij} \|_{L^q_w} m_{i}^{\frac{\mu}{2} \left[ 1- N\left( \frac12-\frac{(2q-1)}{2q \mu}\right) \right]}   m_{j}^{\frac{\mu}{2} \left[ 1- N\left( \frac12-\frac{(2q-1)}{2q \mu}\right) \right]} \  \left\|  \Phi(t) \right\|_{\mathcal{H}^1_A}^{2N \mu \left( \frac12-\frac{(2q-1)}{2q \mu} \right)}  \\
 \noalign{\vskip4pt}
& \leq  C_3 \left\|  \Phi(t) \right\|_{\mathcal{H}^1_A}^{2N \mu \left( \frac12-\frac{(2q-1)}{2q \mu} \right)}.
\end{align*}
Consequently,
 \begin{equation*}
\| \Phi(t) \|^2_{\mathcal{H}^1_A}  \leq C_0 + C_1+ C_2(\varepsilon)+ K \left(  \sum_{j=1}^{m} \frac{N l_j}{4} \varepsilon^{\frac{4}{N l_j}} \right)    \left\|  \Phi(t) \right\|_{\mathcal{H}^1_A}^2 +  C_3 \left\|  \Phi(t) \right\|_{\mathcal{H}^1_A}^{2N \mu \left( \frac12-\frac{(2q-1)}{2q \mu} \right)}
\end{equation*}
and taking $\varepsilon$ such that $1-  K \left(  \sum_{j=1}^{m} \frac{N l_j}{4} \varepsilon^{\frac{4}{N l_j}} \right)  $ is positive, we have
 \begin{equation*}
\Big( 1-  K \Big(  \sum_{j=1}^{m} \frac{N l_j}{4} \varepsilon^{\frac{4}{N l_j}} \Big) \Big)  \left\|  \Phi(t) \right\|_{\mathcal{H}^1_A}^2 \leq C_0 + C_1+ C_2(\varepsilon)+ C_3 \left\|  \Phi(t) \right\|_{\mathcal{H}^1_A}^{2N \mu \left( \frac12-\frac{(2q-1)}{2q \mu} \right)}.
\end{equation*}
By the hypothesis on $\mu$, if  $\left\|  \Phi(t) \right\|_{\mathcal{H}^1_A}^2$ was unbounded respect to $t$,  by the above inequality we would have a contradiction. So we have proved  that $ \left\|  \Phi(t) \right\|_{\mathcal{H}^1_A}^2 $ is bounded respect to $t$ thus proving the global existence result \ref{solGlobal} by standard arguments.
\end{proof}

\begin{remark}
We observe in particular that condition $l_j  < \frac{4}{N}$  in (G) and  $\mu < 2-\frac{1}{q}+\frac{2}{N}$ are fundamental for the proof of the above global existence result.
\end{remark}










\section{Appendix}

\noindent \textbf{Step 1: Construction of approximate solutions (Cazenave-Weissler \cite{CazWeissler}).} We present a supposed adaptation of the arguments in Cazenave and Weissler \cite{CazWeissler} (see Theorem 2.1) in the case of systems (see Remark 2.7 in \cite{CazWeissler}) and  for $A \neq 0$ and the arguments in Cazenave-Esteban \cite{CazEsteban} in the case of systems and not necessarily  for magnetic potentials $A$ of polynomial type and constant magnetic fields $B$. \\
We apply Lemma \ref{embeddH1A}, i.e., $\mathcal{H}^1_A \subset \mathcal{L}^p$ and $\mathcal{L}^{p'} \subset \mathcal{H}'_A$ in the place of usual Sobolev's inequalities and we use suitable estimates of the propagator $T(t)=e^{-itL_A}$. \\
By assumptions (V), (W) and (g) and (h), we have that each $\tilde{g}_{i}$ and so $\tilde{g}$ belong to $C(\mathcal{H}^1_A , \mathcal{H}'_A )$.
Consider the problem
\begin{equation}
    \label{Approxprobl}
\begin{cases}
\iu \de_t \Phi^n - L_A \Phi^n + \tilde{g}_n(\Phi^n)=0, & \\
\noalign{\vskip4pt}
\Phi^n(0)=\varphi \in \mathcal{H}_A^1
\end{cases}
 \end{equation}
 where $\tilde{g}_n= \sum_{k=1}^3 \tilde{g}_{k,n}$ and  each of the $\tilde{g}_{k,n}$ is the natural regularization for every  given type of nonlinearity we deal with as Examples 1, 2 and 3 in \cite{CazWeissler}. So each $\tilde{g}_{k,n}$, $\tilde{G}_{k,n}$ and, by Example 4 in \cite{CazWeissler},  $\tilde{g}_n$ and $\tilde{G}_n$ satisfy all the assumptions in Section 2 of \cite{CazWeissler}. So we are able to prove the following

\begin{lemma}\label{lemmataumA}
  Let $M \geq 0$ and $\varphi \in \mathcal{H}_A^1$ such that $\| \varphi   \|_{\mathcal{H}_A^1} \leq M$. Then,  there exists $\tau_{n,A} > 0  $ such that there exist a sequence $(\Phi^n)_{n \in \N} $ of functions of $C( [0, \tau_{n,A}[,  \mathcal{H}_A^1  ) \cap C^1 ( [0, \tau_{n,A}[, \mathcal{H}'_A  )$ solutions of (\ref{Approxprobl}). Furthermore,  for any $t \in [0, \tau_{n,A}[$, we have
 \begin{equation}\label{conservFAm}
F_{A, n}(\Phi^n)=F_{A, n} (\varphi)
 \end{equation}
 \begin{equation}\label{conservL2}
 \|\Phi^n(t)  \|_{\mathcal{L}^2}= \|\varphi \|_{\mathcal{L}^2}
 \end{equation}
 \end{lemma}

 \begin{proof}
We expect that the proof is the same as in Lemma 2.7 in \cite{LMichael}  in the case of  systems which is  obtained by Lemma 3.5  in \cite{CazWeissler}  replacing usual derivatives by magnetic ones. In particular, since $\tilde{g}_n$ is a globally Lipschitz- continuous nonlinearity,  we can apply the  classical result on $T(t)$ that generates the solution $\Phi^n$ above,  contained in the  Appendix in \cite{CazWeissler} that we recall in the following.
 \medskip

\noindent \textbf{Remark:}
 Let $X$ be a Banach space and $A$ the generator of a $C_0$ semigroup $T(t)$. Let $F \in C(X,X)$ be Lipschitz continuous on bounded sets of $X$. It is well-known that for any $\varphi \in X$, there exists $T_{\max}(\varphi) > 0$ and a unique solution $u \in C([0, T_{\max}], X)$ of

\begin{equation}
u(t)= T(t) \varphi + \int_0^t T(t-s) F(u(s)) \, ds, \quad\hbox{for $t \in [0, T_{\max})$.}
\end{equation}
Moreover, if $T < \infty$ then $\| u(t) \|_{\max} \rightarrow + \infty$ as $t \uparrow  T_{\max}$. Furthermore, the mapping $\varphi \rightarrow T_{\max} (\varphi)$ is lower semicontinuous. If $T \in (0, T_{\max}(\varphi)$ and if $\varphi^n \rightarrow \varphi$ in $X$ as $n \rightarrow \infty$, then $u_{\varphi^n} \rightarrow u_{\varphi}$ in $C([0, T], X)$. If $X$ is reflexive and $\varphi \in D(A)$, then  $u_{\varphi} \in C[0, T), D(A)) \cap
C([0, T), X)$ and $u_{\varphi}$ solves the problem
$$u_t= Au+F(u) \quad\hbox{for $t \in [0, T_{\max})$ and $u(0)=\varphi$}.   $$

\medskip

 \end{proof}

 \bigskip

\noindent \textbf{Step 2: Boundedness of the existence time.} From the conservation laws  (\ref{conservFAm}) and (\ref{conservL2}) of the approximate problem (\ref{Approxprobl}), we show that the existence time $\tau_{n,A}$ can be bounded from below uniformly with respect to $n \in \N$ and $A$ satisfying the Assumptions (A) and (B).

\begin{lemma}\label{lemmaNormaInf}
Let $M > 0$ and let $A$ satisfy assumptions as above with some constants $(C_{\alpha})_{\alpha \in \N^N}$. Then, there exists $T(M) > 0$ depending only on $M$ and the $(C_{\alpha})$'s such that for all $\varphi \in \mathcal{H}_A^1$  such that $\|\varphi \|_{\mathcal{H}_A^1} \leq M$ we have
\begin{equation}
\| \Phi^n \|_{L^{\infty} ([0, T(M)]; \mathcal{H}_A^1)  } \leq 2M= 2\|\varphi \|_{\mathcal{H}_A^1}
\end{equation}
\begin{equation}
\| \Phi^n_t \|_{L^{\infty} ([0, T(M)]; \mathcal{H}'_A)  } \leq C(M)
\end{equation}
\end{lemma}

 \begin{proof}
We expect that the proof is exactly the same as in Lemma 2.8 in \cite{LMichael} in the case of  systems by using Lemma  \ref{lemmataumA} (in particular, we use strongly the conservation of energy (\ref{conservFAm})). Recall  that, in Lemma \ref{embeddH1A}, the constant $C$ is independent on $A$ and  by  a result like Lemma 3.3 in \cite{CazWeissler}  we get uniformity with respect to $A$.
 \end{proof}

\bigskip
\bigskip

\noindent \textbf{Step 3: Passage to the limit.} The final step is to prove the convergence of the $\Phi^n$ to a solution of the initial problem. First, we prove convergence in $\mathcal{L}^2$.

\begin{lemma}
Let $M > 0$ and let $A$ satisfy assumptions as above with some constants $(C_{\alpha})_{\alpha \in \N^N}$. Then, there exists $T(M) > 0$ depending only on $M$ and the $(C_{\alpha})$'s such that for all $\varphi \in \mathcal{H}_A^1$  such that $\|\varphi \|_{\mathcal{H}_A^1} \leq M$, $(\Phi^n)_{n \in \N^N}$ is a Cauchy sequence in $C((0, T(M)]; \mathcal{L}^2 )$.

\end{lemma}

\begin{proof}
We expect that the proof is the same as in Lemma 2.9 in \cite{LMichael} in the case of  systems, using Theorem 2, Lemma 2.4, Lemma 2.3 in \cite{LMichael} or in general Lemma 3.3 in \cite{CazWeissler}   and Lemma \ref{lemmataumA}.
\end{proof}

\bigskip

\noindent We complete the proof of Theorem \ref{solCC1}. We denote by $\Phi$ the limit of $\Phi^n$ in  $C((0, T(M)]; \mathcal{L}^2 )$. From Lemma \ref{lemmaNormaInf}, it follows that $\Phi \in L^{\infty}((-T(M), T(M)); \mathcal{H}^1_A )$ and by Lemma  \ref{embeddH1A}, $\Phi^n$ converges to $\Phi$ in $C((0, T(M)]; \mathcal{L}^r )$ for all $r \in [2, 2^*)$. hence, it follows from a result like Lemma 3.3 in \cite{CazWeissler}  that $\tilde{g}_{n}(\Phi^n)$ converges to $\tilde{g}(\Phi)$ in $C((0, T(M)]; \mathcal{H}'_A )$ and $\Phi$ solves (\ref{Approxprobl}) in $L^{\infty}((0, T(M)]; \mathcal{H}'_A )$.       Uniqueness  is an immediate consequence of Proposition \ref{stimePropagYajima}. Conservation laws are obtained from the passage to the limit. Indeed, combining Lemma 3.3 in \cite{CazWeissler} and Lemma \ref{lemmataumA} we get
$$F_A(\Phi(t))=F_A(\varphi).   $$
This shows that $\Phi \in C((0, T(M)]; \mathcal{H}^1_A )$ and hence  $\Phi \in C^1((0, T(M)]; \mathcal{H}'_A )$. Theorem \ref{solCC1} follows from con\-si\-de\-ring the maximal solution corresponding to the initial datum and the reverse equation.


\bigskip
\bigskip

\end{document}